\documentclass[reqno,11pt]{amsart}
\usepackage{amsmath}
\usepackage{amsxtra}
\usepackage{amscd}
\usepackage{amsthm}
\usepackage{amsfonts}
\usepackage{amssymb}
\usepackage{eucal}
\usepackage[matrix,arrow,curve]{xy}
\theoremstyle{plain}
\newtheorem{Thm}[subsection]{Theorem}
\newtheorem{Cor}[subsection]{Corollary}
\newtheorem{Lem}[subsection]{Lemma}
\newtheorem{Prop}[subsection]{Proposition}
\newtheorem{Conj}[subsection]{Conjecture}

\theoremstyle{definition}
\newtheorem{Def}[subsection]{Definition}

\theoremstyle{remark}

\newtheorem{Rem}[subsection]{Remark}

\newtheorem{conj}{Conjecture}

\numberwithin{equation}{section}
\renewcommand{\rm}{\normalshape}

\newif\ifShowLabels
\ShowLabelstrue
\newdimen\theight
\def\TeXref#1{%
    \leavevmode\vadjust{\setbox0=\hbox{{\tt
        \quad\quad  {\small \rm #1}}}%
    \theight=\ht0
    \advance\theight by \lineskip
    \kern -\theight \vbox to
    \theight{\rightline{\rlap{\box0}}%
    \vss}%
    }}%

\ShowLabelsfalse% comment this out if labels should be printed

\renewcommand{\sec}[2]{\section{#2}\label{S:#1}%
    \ifShowLabels \TeXref{{S:#1}} \fi}
\newcommand{\ssec}[2]{\subsection{#2}\label{SS:#1}%
    \ifShowLabels \TeXref{{SS:#1}} \fi}

\newcommand{\refs}[1]{Section ~\ref{S:#1}}
\newcommand{\refss}[1]{Section ~\ref{SS:#1}}

\newcommand{\reft}[1]{Theorem ~\ref{T:#1}}

\newcommand{\refp}[1]{Proposition ~\ref{P:#1}}
\newcommand{\refc}[1]{Corollary ~\ref{C:#1}}

\newcommand{\refe}[1]{\eqref{E:#1}}
\newcommand{\refco}[1]{Conjecture ~\ref{Co:#1}}

\newenvironment{thm}[1]%
    { \begin{Thm} \label{T:#1}  \ifShowLabels \TeXref{T:#1} \fi }%
    { \end{Thm} }

\renewcommand{\th}[1]{\begin{thm}{#1} \sl }
\renewcommand{\eth}{\end{thm} }

\newenvironment{lemma}[1]%
    { \begin{Lem} \label{L:#1}  \ifShowLabels \TeXref{L:#1} \fi }%
    { \end{Lem} }
\newcommand{\lem}[1]{\begin{lemma}{#1} \sl}
\newcommand{\elem}{\end{lemma}}

\newenvironment{propos}[1]%
    { \begin{Prop} \label{P:#1}  \ifShowLabels \TeXref{P:#1} \fi }%
    { \end{Prop} }
\newcommand{\prop}[1]{\begin{propos}{#1}\sl }
\newcommand{\eprop}{\end{propos}}

\newenvironment{corol}[1]%
    { \begin{Cor} \label{C:#1}  \ifShowLabels \TeXref{C:#1} \fi }%
    { \end{Cor} }
\newcommand{\cor}[1]{\begin{corol}{#1} \sl }
\newcommand{\ecor}{\end{corol}}

\newenvironment{defeni}[1]%
    { \begin{Def} \label{D:#1}  \ifShowLabels \TeXref{D:#1} \fi }%
    { \end{Def} }
\newcommand{\defe}[1]{\begin{defeni}{#1} \sl }
\newcommand{\edefe}{\end{defeni}}

\newenvironment{remark}[1]%
    { \begin{Rem} \label{R:#1}  \ifShowLabels \TeXref{R:#1} \fi }%
    { \end{Rem} }
\newcommand{\rem}[1]{\begin{remark}{#1}}
\newcommand{\erem}{\end{remark}}

\newenvironment{conjec}[1]%
    { \begin{Conj} \label{Co:#1}  \ifShowLabels \TeXref{Co:#1} \fi }%
    { \end{Conj} }
\renewcommand{\conj}[1]{\begin{conjec}{#1} \sl }
\newcommand{\econj}{\end{conjec}}

\newcommand{\eq}[1]%
    { \ifShowLabels \TeXref{E:#1} \fi
       \begin{equation} \label{E:#1} }
\newcommand{\eeq}{ \end{equation} }

\newcommand{\prf}{ \begin{proof} }
\newcommand{\epr}{ \end{proof} }

\newcommand\alp{\alpha}

\newcommand\lam{\lambda}        \newcommand\Lam{\Lambda}

\newcommand\calK{{\mathcal{K}}}

\newcommand\calO{{\mathcal{O}}}

\newcommand\calS{{\mathcal{S}}}
\newcommand\calT{{\mathcal{T}}}

\newcommand\calW{{\mathcal{W}}}

\newcommand\GG{\mathbb{G}}

\newcommand\ZZ{\mathbb{Z}}

\newcommand\CC{\mathbb{C}}

 \newcommand\grg{{\mathfrak{g}}}

\newcommand\sdp{\times \hskip -0.3em {\raise 0.3ex
\hbox{$\scriptscriptstyle |$}}} % semidirect product

\newcommand\codim{\operatorname{codim}}

\newcommand\End{\operatorname{End\,}}
\newcommand\Ext{\operatorname{Ext}}

\newcommand\GL{\operatorname{GL}}
\newcommand\Gr{\operatorname{Gr}}

\newcommand\Hom{\operatorname {Hom}}

\newcommand\id{\operatorname{id}}
\newcommand\Id{\operatorname{Id}}

\newcommand\Int{\operatorname{Int}}

\newcommand\Perv{\operatorname{Perv}}

\newcommand\SL{{\rm SL}}

\newcommand\olam{{\overline{\lambda}}}

\newcommand\hatG{{\widehat{G}}}

\newcommand\tilG{{\widetilde{G}}}

\newcommand\x{\times}

\newcommand{\ra}{\rangle}
\newcommand{\la}{\langle}

\renewcommand{\Id}{\text{Id}}

\newcommand\nc{\newcommand}

\newcommand{\go}{G(\calO)}

\newcommand{\rc}{\rho^{\vee}}
\newcommand{\ogg}{{\overline \gg}}
\newcommand{\ogl}{\ogg^{\lam}}
\newcommand{\pgg}{\text{Perv}_{\go}(\gg)}
\newcommand{\IC}{{\operatorname{IC}}}

\newcommand{\iso}{{\stackrel{\sim}{\longrightarrow}}}

\nc\aff{\operatorname{aff}}
\nc\oGr{\overline{\Gr}}
\nc\Bun{\operatorname{Bun}}
\nc\hgrg{\widehat{\grg}}
\renewcommand\Int{\operatorname{Int}}
\nc\bInt{\overline{\Int}}
\nc\hatLam{\widehat{\Lam}}
\nc\bmu{\overline{\mu}}
\nc\bnu{\overline{\nu}}
\nc\blambda{\overline{\lam}}
\renewcommand\SL{\operatorname{SL}}
\nc\ocalW{\overline{\calW}}
\nc\pos{\operatorname{pos}}
\nc\IH{\operatorname{IH}}
\nc\Rep{\operatorname{Rep}}
\nc\Gal{\operatorname{Gal}}
\nc{\tilGr}{\widetilde{\Gr}}

\nc\Pic{\operatorname{Pic}}
\nc{\HC}{{\mathcal{HC}}}
\nc{\on}{\operatorname}
\nc{\BA}{{\mathbb{A}}}
\nc{\BC}{{\mathbb{C}}}
\nc{\BM}{{\mathbb{M}}}
\nc{\BN}{{\mathbb{N}}}
\nc{\BP}{{\mathbb{P}}}
\nc{\BR}{{\mathbb{R}}}
\nc{\BZ}{{\mathbb{Z}}}
\nc{\BS}{{\mathbb{S}}}

\nc{\CA}{{\mathcal{A}}}
\nc{\CB}{{\mathcal{B}}}
\nc{\CalD}{{\mathcal D}}
\nc{\CE}{{\mathcal{E}}}
\nc{\CF}{{\mathcal{F}}}
\nc{\CG}{{\mathcal{G}}}
\nc{\CH}{{\mathcal{H}}}
\nc{\CK}{{\mathcal{K}}}
\nc{\CL}{{\mathcal{L}}}
\nc{\CM}{{\mathcal{M}}}
\nc{\CMM}{{\mathcal{M}^{\operatorname{gen}}_\hbar(-\rho)}}
\nc{\CN}{{\mathcal{N}}}
\nc{\CO}{{\mathcal{O}}}
\nc{\CP}{{\mathcal{P}}}
\nc{\CQ}{{\mathcal{Q}}}
\nc{\CR}{{\mathcal{R}}}
\nc{\CS}{{\mathcal{S}}}
\nc{\CT}{{\mathcal{T}}}
\nc{\CU}{{\mathcal{U}}}
\nc{\CV}{{\mathcal{V}}}
\nc{\CW}{{\mathcal{W}}}
\nc{\CX}{{\mathcal{X}}}
\nc{\CZ}{{\mathcal{Z}}}

\nc{\gen}{{\operatorname{gen}}}
\nc{\cM}{{\check{\mathcal M}}{}}
\nc{\csM}{{\check{\mathcal A}}{}}
\nc{\obM}{{\overset{\circ}{\mathbf M}}{}}
\nc{\oCA}{{\overset{\circ}{\mathcal A}}{}}
\nc{\obA}{{\overset{\circ}{\mathbf A}}{}}
\nc{\ooM}{{\overset{\circ}{M}}{}}
\nc{\osM}{{\overset{\circ}{\mathsf M}}{}}
\nc{\vM}{{\overset{\bullet}{\mathcal M}}{}}
\nc{\nM}{{\underset{\bullet}{\mathcal M}}{}}
\nc{\obD}{{\overset{\circ}{\mathbf D}}{}}
\nc{\cp}{{\overset{\circ}{\mathbf p}}{}}
\nc{\ofZ}{{\overset{\circ}{\mathfrak Z}}{}}

\nc{\fa}{{\mathfrak{a}}}
\nc{\fb}{{\mathfrak{b}}}
\nc{\fg}{{\mathfrak{g}}}
\nc{\fgl}{{\mathfrak{gl}}}
\nc{\fh}{{\mathfrak{h}}}
\nc{\fj}{{\mathfrak{j}}}
\nc{\fm}{{\mathfrak{m}}}
\nc{\fn}{{\mathfrak{n}}}
\nc{\fu}{{\mathfrak{u}}}
\nc{\fp}{{\mathfrak{p}}}
\nc{\frr}{{\mathfrak{r}}}
\nc{\fs}{{\mathfrak{s}}}
\nc{\fT}{{\mathfrak{T}}}
\nc{\ofT}{{\overline{\mathfrak T}}}
\nc{\ofS}{{\overline{\mathfrak S}}}
\nc{\fsl}{{\mathfrak{sl}}}
\nc{\hsl}{{\widehat{\mathfrak{sl}}}}
\nc{\hgl}{{\widehat{\mathfrak{gl}}}}
\nc{\hg}{{\widehat{\mathfrak{g}}}}
\nc{\chg}{{\widehat{\mathfrak{g}}}{}^\vee}
\nc{\hn}{{\widehat{\mathfrak{n}}}}
\nc{\chn}{{\widehat{\mathfrak{n}}}{}^\vee}

\nc{\fA}{{\mathfrak{A}}}
\nc{\fB}{{\mathfrak{B}}}
\nc{\fD}{{\mathfrak{D}}}
\nc{\fE}{{\mathfrak{E}}}
\nc{\fF}{{\mathfrak{F}}}
\nc{\fG}{{\mathfrak{G}}}
\nc{\fK}{{\mathfrak{K}}}
\nc{\fL}{{\mathfrak{L}}}
\nc{\fM}{{\mathfrak{M}}}
\nc{\fN}{{\mathfrak{N}}}
\nc{\frP}{{\mathfrak{P}}}
\nc{\fS}{{\mathfrak S}}
\nc{\fU}{{\mathfrak{U}}}
\nc{\fZ}{{\mathfrak{Z}}}

\nc{\bb}{{\mathbf{b}}}
\nc{\bc}{{\mathbf{c}}}
\nc{\be}{{\mathbf{e}}}
\nc{\bj}{{\mathbf{j}}}
\nc{\bn}{{\mathbf{n}}}
\nc{\bp}{{\mathbf{p}}}
\nc{\bq}{{\mathbf{q}}}
\nc{\bv}{{\mathbf{v}}}
\nc{\bx}{{\mathbf{x}}}
\nc{\by}{{\mathbf{y}}}
\nc{\bw}{{\mathbf{w}}}
\nc{\bA}{{\mathbf{A}}}
\nc{\bB}{{\mathbf{B}}}
\nc{\bC}{{\mathbf{C}}}
\nc{\bK}{{\mathbf{K}}}
\nc{\bD}{{\mathbf{D}}}
\nc{\bH}{{\mathbf{H}}}
\nc{\bM}{{\mathbf{M}}}
\nc{\bN}{{\mathbf{N}}}
\nc{\bS}{{\mathbf{S}}}
\nc{\bT}{{\mathbf{T}}}
\nc{\bV}{{\mathbf{V}}}
\nc{\bW}{{\mathbf{W}}}
\nc{\bX}{{\mathbf{X}}}
\nc{\bP}{{\mathbf{P}}}
\nc{\bZ}{{\mathbf{Z}}}

\nc{\sA}{{\mathsf{A}}}
\nc{\sB}{{\mathsf{B}}}
\nc{\sC}{{\mathsf{C}}}
\nc{\sD}{{\mathsf{D}}}
\nc{\sF}{{\mathsf{F}}}
\nc{\sK}{{\mathsf{K}}}
\nc{\sM}{{\mathsf{M}}}
\nc{\sO}{{\mathsf{O}}}
\nc{\sQ}{{\mathsf{Q}}}
\nc{\sP}{{\mathsf{P}}}
\nc{\sV}{{\mathsf{V}}}
\nc{\sZ}{{\mathsf{Z}}}
\nc{\sfp}{{\mathsf{p}}}
\nc{\sr}{{\mathsf{r}}}
\nc{\sfb}{{\mathsf{b}}}
\nc{\sfc}{{\mathsf{c}}}
\nc{\sd}{{\mathsf{d}}}
\nc{\sfl}{{\mathsf{l}}}

\nc{\BK}{{\bar{K}}}

\nc{\tA}{{\widetilde{\mathbf{A}}}}
\nc{\tB}{{\widetilde{\mathcal{B}}}}
\nc{\tg}{{\widetilde{\mathfrak{g}}}}
\nc{\tG}{{\widetilde{G}}}
\nc{\TM}{{\widetilde{\mathbb{M}}}{}}
\nc{\tO}{{\widetilde{\mathsf{O}}}{}}
\nc{\tU}{{\widetilde{\mathfrak{U}}}{}}
\nc{\TZ}{{\tilde{Z}}}
\nc{\tx}{{\tilde{x}}}
\nc{\tbv}{{\tilde{\bv}}}
\nc{\tfP}{{\widetilde{\mathfrak{P}}}{}}
\nc{\tz}{{\tilde{\zeta}}}
\nc{\tmu}{{\tilde{\mu}}}

\nc{\td}{\ddot{\underline{d}}{}}
\nc{\tzeta}{\widetilde{\zeta}{}}
\nc{\hd}{{\widehat{\underline{d}}}}
\nc{\hG}{{\widehat{G}}}
\nc{\hBP}{\widehat{\mathbb P}{}}
\nc{\hQ}{{\widehat{Q}}}
\nc{\hsM}{\widehat{\mathsf M}{}}
\nc{\hfM}{\widehat{\mathfrak M}{}}
\nc{\hCP}{\widehat{\mathcal P}{}}
\nc{\hCR}{\widehat{\mathcal R}{}}
\nc{\hCS}{{\widehat{\mathcal S}}}
\nc{\hfZ}{\widehat{\mathfrak Z}{}}

\nc{\urho}{\underline{\rho}}
\nc{\uB}{\underline{B}}
\nc{\uC}{{\underline{\mathbb{C}}}}
\nc{\ui}{\underline{i}}
\nc{\ofP}{{\overline{\mathfrak{P}}}}
\nc{\hrho}{{\hat{\rho}}}

\nc{\unl}{\underline}
\nc{\ol}{\overline}
\nc{\one}{{\mathbf{1}}}
\nc{\two}{{\mathbf{t}}}

\nc{\Tot}{{\mathop{\operatorname{\rm Tot}}}}
\nc{\Hilb}{{\mathop{\operatorname{\rm Hilb}}}}
\nc{\CHom}{{\mathop{\operatorname{{\mathcal{H}}\it om}}}}
\nc{\defi}{{\mathop{\operatorname{\rm def}}}}
\nc{\length}{{\mathop{\operatorname{\rm length}}}}
\nc{\Cliff}{{\mathsf{Cliff}}}
\nc{\Fl}{{\mathsf{Fl}}}
\nc{\Fib}{{\mathsf{Fib}}}
\nc{\Coh}{{\mathsf{Coh}}}
\nc{\FCoh}{{\mathsf{FCoh}}}

\nc{\reg}{{\text{\rm reg}}}

\nc{\cplus}{{\mathbf{C}_+}}
\nc{\cminus}{{\mathbf{C}_-}}
\nc{\cthree}{{\mathbf{C}_*}}
\nc{\Qbar}{{\bar{Q}}}

\nc{\bh}{{\bar{h}}}
\nc{\bOmega}{{\overline{\Omega}}}
\nc\tGr{\widetilde{\Gr}}

\nc{\seq}[1]{\stackrel{#1}{\sim}}

\renewcommand\gg{\Gr_G}
\nc\uS{\underline{S}}
\begin{document}
\title{Pursuing the double affine Grassmannian III: Convolution with affine Zastava}
\author{Alexander Braverman and Michael Finkelberg}% and Alexander Kuznetsov}

\begin{abstract}This is the third paper of a series (started
by~\cite{BF},~\cite{BF2}) which describes
a conjectural analog of the affine Grassmannian for affine Kac-Moody groups
(also known as the double affine Grassmannian).
The current paper is dedicated to describing a conjectural analog of the
convolution diagram for the double affine Grassmannian and affine Zastava.
\end{abstract}
\maketitle
\sec{int}{Introduction}

\ssec{21}{The usual affine Grassmannian}Let $G$ be a connected complex
reductive group with a Cartan torus $T$, and
let $\calK=\CC((s))$, $\calO=\CC[[s]]$. By the {\it affine
Grassmannian} of $G$ we shall mean the quotient
$\gg=G(\calK)/G(\calO)$. It is known (cf. \cite{BD,MV}) that
$\gg$ is the set of $\CC$-points of an ind-scheme over
$\CC$, which we will denote by the same symbol. Note that $\gg$ is defined
for any (not necessarily reductive) group $G$.

Let $\Lam=\Lam_G$ denote the coweight lattice of $G$ and let $\Lam^{\vee}$
denote the dual lattice (this is the weight lattice of $G$).
We let $2\rho_G^{\vee}$ denote the sum of the positive roots of $G$.

The group-scheme $G(\calO)$ acts on $\gg$ on the left and
its orbits can be described as follows.
One can identify the lattice $\Lam_G$ with
the quotient $T(\calK)/T(\calO)$. Fix $\lam\in\Lam_G$ and
let $s^{\lam}$ denote any lift of $\lam$ to $T(\calK)$.
Let $\gg^{\lam}$ denote the $\go$-orbit of $s^{\lam}$
(this is clearly independent of the choice of $s^{\lam}$).
The following result is well-known:
\lem{gras-orbits}
\begin{enumerate}
\item
$$
\gg=\bigcup\limits_{\lam\in\Lam_G}\gg^{\lam}.
$$
\item
We have $\Gr_G^{\lam}=\Gr_G^{\mu}$ if an only if $\lam$ and $\mu$ belong
to the same $W$-orbit on $\Lam_G$ (here $W$ is the Weyl group of $G$).
In particular,
$$
\gg=\bigsqcup\limits_{\lam\in\Lam^+_G}\gg^{\lam}.
$$
\item
For every $\lam\in\Lam^+$ the orbit
$\gg^{\lam}$ is finite-dimensional and its dimension is
equal to $\la\lam,2\rho_G^{\vee}\ra$.
\end{enumerate}
\elem
Let $\ogl$ denote the closure of $\gg^{\lam}$ in $\gg$;
this is an irreducible projective algebraic variety; one has
$\gg^{\mu}\subset \ogl$ if and only if $\lam-\mu$ is a sum of positive roots of
the Langlands dual group $G^{\vee}$.
We will denote by $\IC^{\lambda}$ the intersection
cohomology complex on $\ogl$. Let $\pgg$ denote the category of
$G(\calO)$-equivariant perverse sheaves on $\gg$. It is known
that every object of this category is a direct sum of the
$\IC^{\lam}$'s.

%----------------------------------------------------------------------------

%------------------------------------------------------------------------------------------------------
%--------------------------------------------------------------------------------------------
\ssec{trans-finite}{Transversal slices}Consider the group $G[s^{-1}]\subset G((s))$; let us denote by
$G[s^{-1}]_1$ the kernel of the natural (``evaluation at $\infty$") homomorphism
$G[s^{-1}]\to G$. For any $\lam\in\Lam$ let $\Gr_{G,\lam}=G[s^{-1}]\cdot s^{\lam}$. Then
it is easy to see that one has
$$
\gg=\bigsqcup\limits_{\lam\in\Lam^+}\Gr_{G,\lam}
$$

Let also $\calW_{G,\lam}$ denote the $G[s^{-1}]_1$-orbit of $s^{\lam}$.
For any $\lam,\mu\in\Lam^+$, $\lam\geq \mu$ set
$$
\Gr^{\lam}_{G,\mu}=\gg^{\lam}\cap \Gr_{G,\mu},\quad
\oGr^{\lam}_{G,\mu}=\oGr_G^{\lam}\cap \Gr_{G,\mu}
$$
and
$$
\calW^{\lam}_{G,\mu}=\gg^{\lam}\cap \calW_{G,\mu},\quad
\ocalW^{\lam}_{G,\mu}=\oGr_G^{\lam}\cap \calW_{G,\mu}.
$$
Note that $\ocalW^{\lam}_{G,\mu}$ contains the point $s^{\mu}$ in it. The variety $\ocalW^{\lam}_{G,\mu}$ can be thought of
as a transversal slice to $\gg^{\mu}$ inside $\oGr_G^{\lam}$ at the point $s^{\mu}$ (cf. \cite{BF}, Lemma 2.9).
%-------------------------------------------------------------------------------------------------------------------------------
\ssec{convolution}{The convolution}
We can regard $G(\calK)$ as a total space of a $G(\calO)$-torsor over $\gg$. In particular,
by viewing another copy of $\gg$ as a $G(\calO)$-scheme, we can form the associated fibration
$$
\gg\star\gg:=G(\calK)\underset{\go}\x\gg=G(\calK)\underset{\go}\x G(\calK)/\go.
$$

One has the natural maps $p,m:\gg\star\gg\to\gg$ defined as follows. Let $g\in G(\calK),x\in\gg$. Then
$$
p(g\times x)=g\, \text{mod}\,\go; \quad m(g\times x)=g\cdot x.
$$

For any $\lam_1,\lam_2\in\Lam_G^+$ let us set $\gg^{\lam_1}\star\gg^{\lam_2}$ to be the corresponding subscheme
of $\gg\star\gg$; this is a fibration over $\gg^{\lam_1}$ with the typical fiber $\gg^{\lam_2}$. Its closure is
$\oGr^{\lam_1}\star\oGr^{\lam_2}$. In addition, we define
$$
(\gg^{\lam_1}\star\gg^{\lam_2})^{\lam_3}=m^{-1}(\gg^{\lam_3})\cap
(\gg^{\lam_1}\star\gg^{\lam_2}).
$$
It is known (cf. \cite{Lu-qan}) that
\eq{Lusestimate}
\dim((\gg^{\lam_1}\star\gg^{\lam_2})^{\lam_3})=
\la \lam_1+\lam_2+\lam_3,\rc_G\ra.
\end{equation}
(It is easy to see that although $\rc_G\in \frac{1}{2}\Lam_G^{\vee}$,
the RHS of \refe{Lusestimate} is an integer whenever the above intersection
is non-empty.)

Starting from any perverse sheaf $\calT$ on $\gg$ and a $G(\calO)$-equivariant perverse
sheaf $\calS$ on $\gg$, we can form their twisted external product $\calT\widetilde{\boxtimes}\calS$ (see e.g.~Section~4 of~\cite{MV}), which will
be a perverse sheaf on $\gg\star\gg$. For two objects $\calS_1,\calS_2\in\pgg$ we define their convolution
$$
\calS_1\star\calS_2=m_!(\calS_1\widetilde{\boxtimes}\calS_2).
$$

The following theorem, which is a categorical version of the Satake equivalence,
is a starting point for this paper, cf. \cite{Lu-qan},\cite{Gi} and \cite{MV}. The best reference
so far is \cite{BD}, Sect. 5.3.

\th{grassmannian}
\begin{enumerate}
\item
Let $\calS_1,\calS_2\in \Perv_{\go}(\gg)$. Then
$\calS_1\star\calS_2\in \Perv_{\go}(\gg)$.
\item
The convolution $\star$ extends to a structure of
a tensor category on $\Perv_{\go}(\gg)$.
\item
As a tensor category, $\Perv_{G(\calO)}(\Gr_G)$ is
equivalent to the category $\text{Rep}(G^{\vee})$. Under this
equivalence, the object $\IC^{\lam}$ goes over to the irreducible representation
$L(\lam)$ of $G^{\vee}$ with highest weight $\lam$).
\end{enumerate}
\eth

\ssec{rev1}{}
%This paper is a continuation of~\cite{BF},~\cite{BF2}. The setup is described
%in~Sections~1.1--1.8 of~\cite{BF2}.
%Recall that as a tensor category, $\Perv_{G(\calO)}(\Gr_G)$ is
%equivalent to the category $\text{Rep}(G^{\vee})$.
The equivalence $\Perv_{G(\calO)}(\Gr_G)\iso\text{Rep}(G^{\vee})$ is given
by a fiber functor~\cite{BD},~\cite{MV} of integration over semiinfinite
orbits. Namely, let $N_-\subset G$ be the unipotent radical of the negative Borel subgroup, and let $\fT_\lambda\subset\Gr_G$ be the orbit of $N_-(\CK)$ through the point $s^\lambda\in\Gr_G$. Then the weight $\lambda$ component of the fiber functor is given by the cohomology with supports in $\fT_\lambda$. Let us recall an equivalent construction of this fiber functor.

From now on we assume that $G$ is almost simple and simply connected. 
We consider a smooth curve $\bC$ of genus 0 with two marked points $0,\infty$. 
Let $\Bun_G$ (resp. $\Bun_B$) stand for the moduli stack of $G$-bundles
(resp. $B$-bundles) on $\bC$. Here $B$ is the positive Borel subgroup of $G$.
The natural morphism $\Bun_B\to\Bun_G$ is {\em not} proper, and Drinfeld has
discovered a natural relative compactification $\overline\Bun_B$ of $\Bun_B$.
It is the moduli stack of the following data:

(a) A $G$-bundle $\CF_G$ on $\bC$;

(b) For each dominant weight $\check\lambda$ of $G$, an invertible subsheaf
$\CL^{\check\lambda}\subset\CV_{\CF_G}^{\check\lambda}$. Here $V^{\check\lambda}$
stands for the irreducible $G$-module with highest weight $\check\lambda$,
and $\CV_{\CF_G}^{\check\lambda}$ stands for the associated vector bundle on
$\bC$.

The collection of invertible subsheaves
$\CL^{\check\lambda}\subset\CV_{\CF_G}^{\check\lambda}$ should satisfy the
{\em Pl\"ucker relations}, that is, for any dominant weights $\check\lambda$
and $\check\mu$, the tensor product $\CL^{\check\lambda}\otimes\CL^{\check\mu}
\subset\CV_{\CF_G}^{\check\lambda}\otimes\CV_{\CF_G}^{\check\mu}$ should coincide
with $\CL^{\check\lambda+\check\mu}\subset\CV_{\CF_G}^{\check\lambda+\check\mu}$
under the natural direct summand embedding $\CV_{\CF_G}^{\check\lambda+\check\mu}
\hookrightarrow\CV_{\CF_G}^{\check\lambda}\otimes\CV_{\CF_G}^{\check\mu}$.

The connected components of $\overline\Bun_B$ are numbered by the coweights
$\lambda\in\Lambda$: for $(\CL^{\check\lambda})\in\overline\Bun_B{}^\lambda$
we have $\deg\CL^{\check\lambda}=-\langle\lambda,\check\lambda\rangle$.

We will denote by $\Lambda^{\on{pos}}\subset\Lambda=\Lambda_G$ the cone of
nonnegative linear combinations of positive coroots of $G$. For every
$\alpha\in\Lambda^{\on{pos}}$ we consider the closed embedding
$i_\alpha:\ \overline\Bun_B\hookrightarrow\overline\Bun_B$
given by sending $(\CF_G,\CL^{\check\lambda}\subset\CV_{\CF_G}^{\check\lambda})$
to $(\CF_G,\CL^{\check\lambda}(-\langle\alpha,\check\lambda\rangle\cdot0)\subset
\CV_{\CF_G}^{\check\lambda})$.

Now let $\overline\CH{}^\lambda_0\to\Bun_G\times\Bun_G$ stand for the
Hecke correspondence at the point $0\in\bC$: the pairs of $G$-bundles
$(\CF_G,\CF'_G)$ together with an isomorphism $\sigma:\ \CF_G\to\CF'_G$ 
off $0\in\bC$ whose pole at $0\in\bC$ has order
less than or equal to $\lambda$. The fibers of the projection $p_1$ (resp.
$p_2$) of $\overline\CH{}^\lambda_0$ to the first (resp. second) copy of
$\Bun_G$ are both isomorphic to $\ogl$.

We define the Hecke correspondence
$(p,\phi):\ \overline\CG{}^\lambda_0=\overline\CH{}^\lambda_0\underset{\Bun_G}
\times\overline\Bun_B\to\overline\Bun_B\times\overline\Bun_B$.
It is the moduli stack of the following data:

(a) a pair of $G$-bundles $\CF_G$ and $\CF'_G$ together with an isomorphism
off $0\in\bC$ lying in $\overline\CH{}^\lambda_0$;

(b) For each dominant weight $\check\lambda$ of $G$, an invertible subsheaf
$\CL^{\check\lambda}\subset\CV_{\CF'_G}^{\check\lambda}$ satisfying the Pl\"ucker
relations.

Forgetting the datum of $\CF_G$ defines the morphism
$p:\ \overline\CG{}^\lambda_0\to\overline\Bun_B$. The morphism
$\phi:\ \overline\CG{}^\lambda_0\to\overline\Bun_B$ is defined as follows.
The condition $(\CF_G,\CF'_G)\in\overline\CH{}^\lambda_0$ implies
$\CV^{\check\lambda}_{\CF'_G}\subset\CV^{\check\lambda}_{\CF_G}(\langle-w_0\lambda,
\check\lambda\rangle\cdot0)$ for every dominant weight $\check\lambda$.
Hence $\CL^{\check\lambda}(\langle w_0\lambda,\check\lambda\rangle\cdot0)\subset
\CV^{\check\lambda}_{\CF_G}$, and we set $\phi(\CF_G,\CF'_G,
\CL^{\check\lambda}\subset\CV_{\CF'_G}^{\check\lambda}):=
(\CF_G,\CL^{\check\lambda}(\langle w_0\lambda,\check\lambda\rangle\cdot0)\subset
\CV^{\check\lambda}_{\CF_G})$.

Finally, we are able to state a theorem (see~\cite{BG}~Theorem~3.1.4
and~\cite{FM}~Theorem~13.2) providing a version of the fiber functor from
the category $\Perv_{G(\calO)}(\Gr_G)$ to $\text{Rep}(G^{\vee})$.
For a finite-dimensional $G^\vee$-module $V$, and $\mu\in\Lambda_G$, we
denote the $\mu$-weight subspace of $V$ by $V(\mu)$.

\th{bgfm}
$\phi_!\IC(\overline\CG{}^\lambda_0)\simeq\bigoplus_{\alpha\in\Lambda_G^{\on{pos}}}
i_{\alpha!}\IC(\overline\Bun_B)\otimes V^\lambda(w_0(\lambda)+\alpha)$.
\eth

\ssec{rev2}{}
The goal of the present paper is to formulate an analogue of~\reft{bgfm}
for the double affine Grassmannian. However, as we have seen
in~\cite{BF},~\cite{BF2}, the affine versions of the objects like $\ogl$ or
$\overline\Bun_B$ are
out of reach at the moment being ``too global'', and have to be replaced
by certain transversal slices.

A transversal slice to the closed embedding
$i_\alpha:\ \overline\Bun_B\hookrightarrow\overline\Bun_B$ is a well known
{\em Drinfeld Zastava} space $Z^\alpha$ (see~\cite{FM},~\cite{BFGM}).
It is defined as the moduli scheme of collections of invertible subsheaves
$\CL^{\check\lambda}\subset\CO_\bC\otimes V^{\check\lambda}$ satisfying
the Pl\"ucker relations, the degree conditions
$\deg\CL^{\check\lambda}=-\langle\alpha,\check\lambda\rangle$, and the
conditions at $\infty\in\bC$: each
$\CL^{\check\lambda}\subset\CO_\bC\otimes V^{\check\lambda}$ is a {\em
line subbundle} near $\infty\in\bC$, and the fiber $\CL^{\check\lambda}|_\infty$
coincides with the highest line in $V^{\check\lambda}$.

By construction, we have a locally closed embedding
$z_\alpha:\ Z^\alpha\hookrightarrow\overline\Bun_B$, and we define the scheme
$\overline{\CG Z}{}^{\lambda,\alpha}$ as the cartesian product of
$z_\alpha:\ Z^\alpha\hookrightarrow\overline\Bun_B$ and
$\phi:\ \overline\CG{}^\lambda_0\to\overline\Bun_B$.

For any $\beta\leq\alpha\in\Lambda^{\on{pos}}_G$ we also have a closed
embedding $i_\beta^\alpha:\ Z^{\alpha-\beta}\hookrightarrow Z^\alpha$ which
sends a collection
$(\CL^{\check\lambda}\subset\CO_\bC\otimes V^{\check\lambda})$ to a
collection $(\CL^{\check\lambda}(-\langle\beta,\check\lambda\rangle\cdot0)
\subset\CO_\bC\otimes V^{\check\lambda})$.

Now~\reft{bgfm} can be equivalently formulated as follows (see~Theorem~13.2
of~\cite{FM}):

\th{fm}
$\phi_!\IC(\overline{\CG Z}{}^{\lambda,\alpha})\simeq
\bigoplus_{\beta\leq\alpha}
i_{\beta!}^\alpha\IC(Z^{\alpha-\beta})\otimes V^\lambda(\lambda-\beta)$.
\eth

The key observation underlying the proof of the theorem is that $\phi^{-1}(i^\alpha_\alpha(0))\cong\ol\fT_{\lambda-\alpha}\cap\ol\Gr{}_G^\lambda$.

\ssec{gaff}{The group $G_{\aff}$}
%From now on we assume that $G$ is almost simple and simply connected.
To a connected reductive group $G$ as above one can associate the corresponding
affine Kac-Moody group $G_{\aff}$ in the following way.
One can consider the polynomial loop group $G[t,t^{-1}]$ (this is an infinite-dimensional group ind-scheme)

It is well-known that $G[t,t^{-1}]$ possesses a canonical central extension $\tilG$ of $G[t,t^{-1}]$:
$$
1\to \GG_m\to \tilG\to G[t,t^{-1}]\to 1.
$$
Moreover,\ $\tilG$ has again a natural structure of a group ind-scheme.

The multiplicative group $\GG_m$ acts naturally on $G[t,t^{-1}]$ and this action lifts to $\tilG$.
We denote the corresponding semi-direct product by $G_{\aff}$; we also let $\grg_{\aff}$ denote its Lie algebra.

The Lie algebra $\grg_{\aff}$ is an untwisted affine Kac-Moody Lie algebra.
In particular,
it can be described by the corresponding affine root system. We denote by $\grg_{\aff}^{\vee}$ the
{\em Langlands dual affine Lie algebra} (which corresponds to the dual affine root system)
and by $G^{\vee}_{\aff}$ the corresponding dual affine Kac-Moody group, normalized by the property
that it contains $G^{\vee}$ as a subgroup (cf. \cite{BF}, Subsection 3.1 for more details).

We denote by $\Lam_{\aff}=\ZZ\x\Lam\x\ZZ$ the coweight lattice of $G_{\aff}$; this is the same as the
weight lattice of $G_{\aff}^{\vee}$. Here the first $\ZZ$-factor
is responsible for the center of $G_{\aff}^{\vee}$ (or $\hatG^{\vee}$);
it can also be thought of as coming from the loop
rotation in $G_{\aff}$. The second $\ZZ$-factor is responsible for the loop rotation in $G_{\aff}^{\vee}$  it may also be thought of
as coming from the center of $G_{\aff}$).
We also denote $\ZZ\times\Lam\subset\ZZ\times\Lam\times\ZZ$ by $\widehat\Lam$,
and we denote $k\times\Lam\subset\ZZ\times\Lam$ by $\widehat\Lam_k\subset
\widehat\Lam$.
We denote by $\Lam_{\aff}^+$ the set of dominant weights of $G_{\aff}^{\vee}$ (which is the same as the set of dominant
coweights of $G_{\aff}$). We also denote by $\Lam_{\aff,k}$ the set of weights of $G_{\aff}^{\vee}$ of level $k$,
i.e. all the weights of the form $(k,\olam,n)$. We put $\Lam_{\aff,k}^+=\Lam_{\aff}^+\cap \Lam_{\aff,k}$.

\smallskip
\noindent
{\bf Important notational convention:}
From now on we shall denote elements
of $\Lam$ by $\blambda,\bmu...$ (instead of just writing $\lam,\mu...$ in order to distinguish them from
the coweights of $G_{\aff}$ (= weights of $G_{\aff}^{\vee}$), which we shall just denote by
$\lam,\mu...$

Let
$\Lam_k^+\subset \Lam$ denote the set of dominant coweights of $G$ such that $\la \blambda,\alp)\leq k$
when $\alp$ is the highest root of $\grg$.
Then it is well-known that a weight $(k,\olam,n)$ of
$G_{\aff}^{\vee}$ lies in $\Lam_{\aff,k}^+$  if and only if $\olam\in\Lam_k^+$ (thus $\Lam_{\aff,k}=\Lam_k^+\x \ZZ$).

Let also $W_{\aff}$ denote affine Weyl group of $G$ which is the semi-direct product of $W$ and $\Lam$.
It acts on the lattice $\Lam_{\aff}$ (resp. $\hatLam$)
preserving each $\Lam_{\aff,k}$ (resp. each $\hatLam_k$). In order to describe this action explicitly it is convenient
to set $W_{\aff,k}=W\ltimes k\Lam$ which naturally acts on $\Lam$. Of course the groups $W_{\aff,k}$ are canonically
isomorphic to $W_{\aff}$ for all $k$. Then the restriction of the $W_{\aff}$-action to $\Lam_{\aff,k}\simeq\Lam\x\ZZ$
comes from the natural $W_{\aff,k}$-action on the first multiple.

It is well known that every $W_{\aff}$-orbit on $\Lam_{\aff,k}$ contains unique element of $\Lam_{\aff,k}^+$.
This is equivalent to saying that $\Lam_k^+\simeq \Lam/W_{\aff,k}$.

\ssec{rev3}{}
%Now we turn to the affine case,
%invoking the notational conventions of~Section~1.7 of~\cite{BF2}.
Our main dream is to create an analog of the affine Grassmannian $\Gr_G$ and the above results
about it in the case when $G$ is replaced by the (infinite-dimensional) group $G_{\aff}$.
The first attempt to do so was made in~\cite{BF}: namely, in {\em loc. cit.} we have constructed analogs of the varieties $\ocalW^{\lam}_{G,\mu}$ in the case
when $G$ is replaced by  $G_{\aff}$. In~\cite{BF2}, we constructed analogs of the varieties
$m_n^{-1}(\ocalW^{\lam}_{G,\mu})\cap (\Gr_G^{\lam_1}\star\cdots\star\Gr_G^{\lam_n})$ and
$m_n^{-1}(\ocalW^{\lam}_{G,\mu})\cap (\oGr_G^{\lam_1}\star\cdots\star\oGr_G^{\lam_n})$ (here
$\lam=\lam_1+\cdots+\lam_n$) when
$G$ is replaced by $G_{\aff}$. We have also constructed analogs of the
corresponding pieces in the {\em Beilinson-Drinfeld Grassmannian} for
$G_{\aff}$.

We will denote by
$\Lambda^{\on{pos}}_{\on{aff}}$ the cone of nonnegative linear combinations
of positive roots of $G^\vee_{\on{aff}}$.
For $\alpha\in\Lambda^{\on{pos}}_{\on{aff}}$ the affine Drinfeld Zastava space
$Z^\alpha$ was constructed in~\cite{BFG}. It is a certain closure of the space
of degree $\alpha$ based maps from $(\bC,\infty)$ to the Kashiwara flag scheme
of $G_{\on{aff}}$. We also have parabolic versions $Z^\theta_{G_{\aff},P}$ of $Z^\alpha=Z^\alpha_{G_{\aff},I}$ ($I$ stands for the Iwahori subgroup of $G_{\aff}$) which are
certain closures of the spaces of based maps from $(\bC,\infty)$ to the
Kashiwara parabolic flag schemes. Among those, the Uhlenbeck space
$\CU_G^a(\BA^2)=Z^a_{G_{\aff},G[t]}$ stands out: it corresponds to the maximal parabolic containing
all the finite simple roots.

Unfortunately, the definition of Zastava given in~\refss{rev2} produces in
the affine case a scheme of {\em infinite type} $\bZ^\alpha$. In the maximal
parabolic case, the Uhlenbeck space $\CU_G^a(\BA^2)$ is a partial resolution
of $\bZ^a_{G_{\aff},G[t]}$. We have a natural forgetting morphism $\bZ^\alpha\to\bZ^a_{G_{\aff},G[t]}$ where
$a$ is the coefficient of the affine simple root in $\alpha$, and
$Z^\alpha$ is defined as the cartesian product of $\bZ^\alpha$ and
$\CU_G^a(\BA^2)$ over $\bZ^a_{G_{\aff},G[t]}$. It is an affine scheme of finite type.

The disadvantage of the above definition is that $Z^\alpha$ does not solve
any moduli problem, and hence is very cumbersome to work with. However, in
the case $G=\SL(N)$, the Zastava space $Z^\alpha$ possesses a semismall
resolution of singularities $\CP^\alpha$,
an {\em affine Laumon space}~\cite{FGK}, which is a moduli space of parabolic
sheaves on $\bC\times\BP^1$. Moreover, according to~\cite{Q}, $\CP^\alpha$
admits a realization as a quiver variety, i.e. as a certain GIT quotient.
The corresponding categorical quotient $\fZ^\alpha$ is an affine reduced
irreducible normal scheme isomorphic to $Z^\alpha$ (see~\cite{Q},\cite{BF4}).
%which is an isomorphism at the level of $\BC$-points (and conjecturally just an
%isomorphism). Since in this paper we are interested in topological questions
%only (such as stalks of IC sheaves), we can replace $Z^\alpha$ by $\fZ^\alpha$
%for all practical purposes.

\ssec{rev4}{}
The main result of the present paper is a construction of an affine version
of the scheme $\overline{\CG Z}{}^{\lambda,\alpha}$ equipped
with a morphism $\phi$ to the affine Zastava $Z^\alpha$.
It is constructed as a quiver variety
in the case $G=\SL(N)$, and then for general $G$ via the adjoint embedding
$G\hookrightarrow\SL(\fg)$. We conjecture that~\reft{fm} holds true in the
affine setting as well.

Although we cannot describe $\overline{\CG Z}{}^{\lambda,\alpha}$ as a
solution of a moduli problem, its open subscheme $\CG Z^{\lambda,\alpha}$
does admit such a description. Let us first assume $\lambda$ has level 1.
We consider the projective plane $\BP^2$ with homogeneous coordinates
$[z_0:z_1:z_2]$ such that the line $\ell_\infty$ ``at infinity'' is given by
the equation $z_0=0$, while $\bC=\ell_0\subset\BP^2$ is given by the
equation $z_2=0$. We consider the blowup $\hBP^2$ at the origin ($z_1=z_2=0$),
and keep the names $\ell_\infty$ and $\ell_0$ for the proper transforms of
$\ell_\infty$ and $\ell_0$. Then $\CG Z^{\lambda,\alpha}$ is the moduli space
of $G$-bundles on $\hBP^2$ equipped with a reduction to $B$ along $\ell_0$
framed at $\ell_\infty$. Note that even with this modular definition, the
construction of projection $\phi:\ \CG Z^{\lambda,\alpha}\to Z^\alpha$ is
rather nontrivial, cf.~\refss{BD}. For an explanation why the moduli space of $G$-bundles
on the blowup appears in the convolution diagram of double affine Grassmannian and affine
Zastava, the interested reader may consult~\cite[Section~8]{BK}.

Let us now assume $\lambda$ has level $k$. We consider the blowup $\hBP^2_k$
of $\BP^2$ at the origin, but not at the maximal ideal of the origin this time;
rather at the ideal generated locally by $(z_1^k,z_2)$. This blowup has an
isolated singularity of Kleinian type $A_{k-1}$ lying off the proper transforms
of $\ell_\infty$ and $\ell_0$. We consider the stacky resolution $\hCS^k$ of
$\hBP^2_k$. Then again $\CG Z^{\lambda,\alpha}$ is the moduli space
of $G$-bundles on $\hCS^k$ equipped with a reduction to $B$ along $\ell_0$
framed at $\ell_\infty$. The projection
$\phi:\ \CG Z^{\lambda,\alpha}\to Z^\alpha$ is constructed in~\refss{pik}.
Similarly to the finite dimensional case of~\refss{rev2}, we have $\phi^{-1}(i^\alpha_\alpha(0))\cong\fT_{\lambda-\alpha}\cap\Gr_{G_{\aff}}^\lambda$; for the definition of the RHS and the proof of the isomorphism, see~\cite[Section~8, especially~Proposition~8.7]{BK}

In the special case when $\lambda-\alpha=\mu:=(k,0,0)$, we have an intermediate open subspace
$\CG Z^{\lambda,\alpha}\subset\ol\CW{}^\lambda_{G_{\aff},\mu}\subset\overline{\CG Z}{}^{\lambda,\alpha}$.
In~Section~3.2 of~\cite{BF3} we have defined the {\em repellents} $\fT^e_\mu\subset\ol\CW{}^\lambda_{G_{\aff},\mu}$; they were also considered in~\cite{Na} under the name of {\em MV cycles}. We conjecture that the central fiber $\phi^{-1}(i^\alpha_\alpha(0))\cap\ol\CW{}^\lambda_{G_{\aff},\mu}$ coincides with $\fT^e_\mu$, and we prove the inclusion $\fT^e_\mu\subset\phi^{-1}(i^\alpha_\alpha(0))\cap\ol\CW{}^\lambda_{G_{\aff},\mu}$ in~\refp{repel}.

\ssec{struc}{Structure of the paper}
In~\refs{zastava} we recall the description of the affine Zastava $Z^\alpha_{\SL(N)_{\aff}}$ in terms of representations of the {\em chainsaw quiver} of~\cite{Q}. Contrary to the ``global" approach of {\em loc. cit.}, we follow the classical ADHM approach on a 2-dimensional toric Deligne-Mumford stack $\CS_N=\BP^1/\unl{\mu}{}_N\times\BP^1$. Here $\unl{\mu}{}_N$ is the group of $N$-th roots of unity, acting on $\BP^1$ with fixed points $0,\infty$, and the quotient is categorical near $\infty$, and stacky near~$0$. In~\refss{gammak}--\refss{def} we describe the irreducible components of the fixed point set $(Z^\alpha_{\SL(N)_{\aff}})^{\Gamma_k}$ of the action of a cyclic group $\Gamma_k$. In the central~\refs{zastava blowup} we describe the parabolic torsion free sheaves on the blowup $\widehat\BP{}^2$ in terms of the {\em dented chainsaw quiver $\widehat Q$}~(\refss{den chain}). The description is modeled on the one in~\cite{NY} for torsion free sheaves on the blowup. The key~\reft{queer blowup} identifying the moduli space of parabolic torsion free sheaves on the blowup with a moduli space of $\widehat Q$-modules is due to A.~Kuznetsov. We introduce the Zastava space for the blowup as the moduli space of $\widehat Q$-modules with certain stability conditions~(\refss{stab}). In~\refs{klein blowup} we introduce the Zastava space for the Kleinian blowup $\widehat\CS{}^k$~(\refss{rev4}) via a trick identifying it with a $\Gamma_k$-fixed points component in the Zastava space on the blowup $\widehat\BP{}^2$. This allows us to describe it as a moduli space of representations of the {\em rift quiver}~(\refss{bicycle},~\reft{maink}). Finally, in~\refs{general G}, for an arbitrary almost simple simply connected group $G$, we define the Zastava space for the Kleinian blowup in terms of the one for $\SL(\fg)$.

\ssec{ack}{Acknowledgments}
It is clear from the above that the paper owes its existence to A.~Kuznetsov's generous 
explanations. We are grateful to him for the permission to reproduce his proof of the 
key~\reft{queer blowup}. As our masters put it, ``Il avait \'et\'e d'abord
pr\'evu que A.~Kuznetsov soit coauteur du pr\'esent article. Il a pr\'ef\'er\'e
s'en abstenir, pour ne pas \^etre corresponsable des erreurs ou
impr\'ecisions qui s'y trouvent. Il n'en est pas moins responsable de bien
des id\'ees que nous exploitons.''
Thanks are due to the referee for the careful reading of the 
manuscript and his valuable comments and suggestions.
We are happy to thank the IAS at the Hebrew University of Jerusalem for the excellent working conditions. 
M.F. was partially supported by the RFBR grants 12-01-00944, 12-01-33101, 
the National Research University Higher School of Economics' Academic Fund 
award No.12-09-0062 and
the AG Laboratory HSE, RF government grant, ag. 11.G34.31.0023.
This study comprises research findings from the ``Representation Theory
in Geometry and in Mathematical Physics" carried out within The
National Research University Higher School of Economics' Academic Fund Program
in 2012, grant No 12-05-0014.

%-------------------------------------------------------------------------------------------------------------------------

\sec{zastava}{Zastava as a quiver variety}

\ssec{chain}{Chainsaw}
We recall some material from~Section~2 of~\cite{Q}. We consider the representations of the following
{\em chainsaw quiver} $Q$

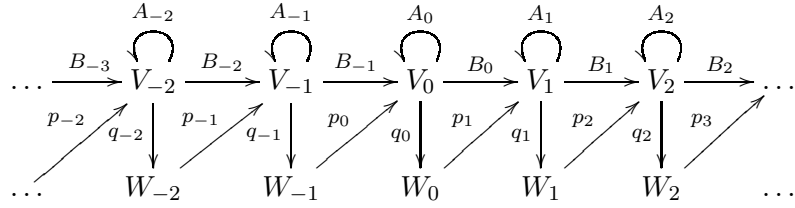
$$\xymatrix{
\ldots \ar[r]^{B_{-3}}
& V_{-2} \ar@(ur,ul)[]_{A_{-2}} \ar[r]^{B_{-2}} \ar[d]_{q_{-2}}
& V_{-1} \ar@(ur,ul)[]_{A_{-1}} \ar[r]^{B_{-1}} \ar[d]_{q_{-1}}
& V_0 \ar@(ur,ul)[]_{A_0} \ar[r]^{B_0} \ar[d]_{q_0}
& V_1 \ar@(ur,ul)[]_{A_1} \ar[r]^{B_1} \ar[d]_{q_1}
& V_2 \ar@(ur,ul)[]_{A_2} \ar[r]^{B_2} \ar[d]_{q_2} &\ldots\\
\ldots \ar[ur]^{p_{-2}} & W_{-2} \ar[ur]^{p_{-1}} & W_{-1} \ar[ur]^{p_0}
& W_0 \ar[ur]^{p_1} & W_1 \ar[ur]^{p_2} & W_2 \ar[ur]^{p_3} &\ldots
}$$
with relations $A_{l+1}B_l-B_lA_l+p_{l+1}q_l=0\ \forall l$. Here the lower
indices run through $\BZ/N\BZ$, and $\dim V_l=d_l,\ \dim W_l=1$.
We denote by $\unl{d}$ the collection of positive integers
$(d_l)_{l\in\BZ/N\BZ}$. We denote by $\sM_{\unl{d}}$ the scheme
of representations of $Q$: a closed subscheme of
$$\bigoplus_{l\in\BZ/N\BZ}\End(V_l)\oplus
\bigoplus_{l\in\BZ/N\BZ}\Hom(V_l,V_{l+1})\oplus
\bigoplus_{l\in\BZ/N\BZ}\Hom(W_{l-1},V_l)\oplus
\bigoplus_{l\in\BZ/N\BZ}\Hom(V_l,W_l)$$
given by equations $A_{l+1}B_l-B_lA_l+p_{l+1}q_l=0\ \forall l$.
We denote by $G_{\unl{d}}$ the group $\prod_{l\in\BZ/N\BZ}\GL(V_l)$;
it acts naturally on $\sM_{\unl{d}}$.
We denote by $\fZ_{\unl{d}}$ the categorical quotient
$\sM_{\unl{d}}/\!/G_{\unl{d}}$. According to~\cite[Theorem~2.7]{Q} and~\cite[Theorem~3.5]{BF4},
$\fZ_{\unl{d}}$ is a reduced irreducible normal scheme isomorphic to
the affine Drinfeld Zastava space $Z^{\unl{d}}_{\SL(N)}$
introduced in~\cite{BFG}.

Furthermore, we consider an open subscheme
$\sM_{\unl{d}}^s\subset\sM_{\unl{d}}$ formed by all the {\em stable}
representations of $Q$, i.e. those $(A_l,B_l,p_l,q_l)_{l\in\BZ/N\BZ}\in
\sM_{\unl{d}}$ such that there is no proper $\BZ/N\BZ$-graded subspace
$V'_\bullet\subset V_\bullet$ stable under $A_\bullet,B_\bullet$ and
containing $p(W_\bullet)$. Then the action of $G_{\unl{d}}$ on
$\sM_{\unl{d}}^s$ is free, and the GIT quotient
$\fM_{\unl{d}}=\sM_{\unl{d}}^s/G_{\unl{d}}$
is a semismall resolution of $\fZ_{\unl{d}}$. Moreover, according to
Section~2.3 of~\cite{Q}, $\fM_{\unl{d}}$ is isomorphic to the moduli space
$\CP_{\unl{d}}$
of torsion free parabolic sheaves of degree $\unl{d}$ on a surface $\bS$.
Here $\bS$ is the product of two projective lines $\bC$ and $\bX$ with
marked points $0_\bX,\infty_\bX\in\bX$ and $0_\bC,\infty_\bC\in\bC$.
The sheaves in question are equipped with a parabolic structure along a line
$\bD_0:=\bC\times0_\bX$, and with a trivialization at ``infinity''
$\bD_\infty:=\bC\times\infty_\bX\cup\infty_\bC\times\bX$.
The isomorphism $\fM_{\unl{d}}\simeq\CP_{\unl{d}}$ is deduced in
{\em loc. cit.} from the ``parabolic vs. orbifold'' correspondence
of~\cite{Bi} by global considerations. We will
rephrase the argument in more local terms in~\refss{local} and~\refss{inverse}
after some preparation in~\refss{adhm} and~\refss{S_N}.

\ssec{adhm}{ADHM}
To warm up we recall the classical ADHM construction
(see e.g.~Section~2 of~\cite{Nak}) following the approach of~Section~5
of~\cite{BGK}.
To this end we introduce the homogeneous
coordinates $(z:t)$ (resp. $(y:x)$) on $\bC$ (resp. $\bX$) such that
$0_\bC$ (resp. $0_\bX$) is given by $z=0$ (resp. $y=0$), and
$\infty_\bC$ (resp. $\infty_\bX$) is given by $t=0$ (resp. $x=0$).
The ADHM construction goes
as follows. We consider the vector spaces $V=\CC^d,\ W=\CC^N$, and the
subscheme $\sM_{N,d}\subset\End(V)\oplus\End(V)\oplus\Hom(W,V)\oplus
\Hom(V,W)$ cut out by the equation $AB-BA+pq=0\ (A,B\in\End(V),\
p\in\Hom(W,V),\ q\in\Hom(V,W))$. We consider an open subscheme
$\sM_{N,d}^s\subset\sM_{N,d}$ formed by all the {\em stable} quadruples
$(A,B,p,q)$, i.e. such that $V$ has no proper subspaces stable under $A,B$
and containing $p(W)$. The group $\GL(V)$ acts naturally on $\sM_{N,d}$;
its action on $\sM_{N,d}^s$ is free, and the GIT quotient
$\sM_{N,d}^s/\GL(V)$ is denoted by $\fM_{N,d}$. It is isomorphic to the
moduli space $\CM_{N,d}$ of torsion free sheaves of rank $N$ and degree $d$
on $\bS$ trivialized at $\bD_\infty$. Namely, $(A,B,p,q)\in\fM_{N,d}$
goes to the middle cohomology of the following monad of vector bundles on
$\bS$:
$$V\otimes\CO_\bS(-1,-1)\stackrel{C}{\longrightarrow}
\begin{array}{c}
V\otimes\CO_\bS(0,-1)\\
\bigoplus\\
V\otimes\CO_\bS(-1,0)\\
\bigoplus\\
W\otimes\CO_\bS
\end{array}
\stackrel{D}{\longrightarrow}
V\otimes\CO_\bS,$$
$C=(tA-z,xB-y,txq),\ D=(-xB+y,tA-z,p)$, where we write $\CO_\bS(-1,-1)$
for $\CO_\bC(-1)\boxtimes\CO_\bX(-1)$, and $\CO_\bS(0,-1)$ for
$\CO_\bC\boxtimes\CO_\bX(-1)$, etc., and we view $x,y$ (resp. $z,t$) as a basis
of $\Gamma(\bX,\CO_\bX(1))$ (resp. $\Gamma(\bC,\CO_\bC(1))$).

\ssec{S_N}{Stack $\CS_N$}
We define a one-dimensional Deligne-Mumford stack $\CX_N$ as follows.
Let $\bX_N\stackrel{\theta}{\to}\bX$ denote the
$N$-fold cyclic covering ramified
over $0_\bX$ and $\infty_\bX$. It is equipped with the action of the
Galois group
$\Gamma_N\simeq\BZ/N\BZ$. The action of $\Gamma_N$ on $\theta^{-1}(\bX-0_\bX-
\infty_\bX)$ is free, and the quotient is $\bX-0_\bX-\infty_\bX$.
We glue the categorical quotient $\theta^{-1}(\bX-0_\bX)/\!/\Gamma_N=\bX-0_\bX$
with the stack quotient $\theta^{-1}(\bX-\infty_\bX)/\Gamma_N$ over the common
open substack $\bX-0_\bX-\infty_\bX$ to obtain the desired stack $\CX_N$.
Note that $\CX_N$ is equipped with a projection $\vartheta$ to $\bX$ which
is an isomorphism off $0_\bX$. The unique point of $\CX_N$ lying over
$0_\bX$ will be denoted by $0_\CX$; its group of automorphisms is $\Gamma_N$.
The unique point of $\CX_N$ lying over $\infty_\bX$ will be denoted by
$\infty_\CX$. Since $N$ is fixed throughout
the Section, we will omit the lower index $N$ to simplify the notations.

We denote $\CO_\CX(\pm\infty_\CX)$ by $\CO_\CX(\pm N)$. For $0\leq l\leq N$
we denote $\CO_\CX(-l\cdot0_\CX)$ by $\CR_l$. Note that
$\CR_{N}\simeq\CO_\CX(-N)$. We have the canonical embeddings
$$\CR_0(-N)\simeq\CR_N\stackrel{\xi_N}{\hookrightarrow}
\CR_{N-1}\stackrel{\xi_{N-1}}{\hookrightarrow}
\ldots\stackrel{\xi_3}{\hookrightarrow}
\CR_2\stackrel{\xi_2}{\hookrightarrow}
\CR_1\stackrel{\xi_1}{\hookrightarrow}\CR_0.$$

We define a 2-dimensional Deligne-Mumford stack $\CS_N$ as
$\bC\times\CX_N$; by an abuse of notation we will denote by $\vartheta$
its projection $\id\times\vartheta$ onto $\bS$.
We denote $\bC\times0_\CX$ by $\CalD_0$, and we denote $\infty_\bC\times\CX
\cup\infty_\CX\times\bC$ by $\CalD_\infty$. By an abuse of notation, we
denote by $\CR_l$ the line bundle $\CO_\bC\boxtimes\CR_l$, and we denote
by $\xi_l$ the morphism $\id\boxtimes\xi_l$.
According to~\cite{Bi}, there
is a one-to-one correspondence between the (torsion free, framed at
$\CalD_\infty$) sheaves on $\CS$, and the (torsion free, framed at
$\bD_\infty$) sheaves on $\bS$ with parabolic structure along $\bD_0$.
Thus $\CP_{\unl{d}}$ is the moduli space of torsion free sheaves of
degree $\unl{d}$ on $\CS$ framed at $\CalD_\infty$.

\ssec{local}{Monad for the stack $\CS_N$}
Finally we are able to recall an ADHM-like construction of the isomorphism
$\fM_{\unl{d}}\iso\CP_{\unl{d}}$. Note that
$\vartheta_*$ establishes an isomorphism $\Gamma(\CX,\CO_\CX(1))\simeq
\Gamma(\bX,\CO_\bX(1))=\BC\langle x,y\rangle$. The desired isomorphism
$\fM_{\unl{d}}\iso\CP_{\unl{d}}$ sends a representative
$(A_\bullet,B_\bullet,p_\bullet,q_\bullet)$ to the middle cohomology of
the following monad of vector bundles on $\CS$:
%(we understand the indices as running through $\BZ/N\BZ$):
\eq{mon}
\bigoplus_{0<l\leq N} V_l\otimes\CR_l\otimes\CO_\CS(-1,0)
\stackrel{C}{\longrightarrow}
\begin{array}{c}
\bigoplus\limits_{0<l\leq N} V_l\otimes\CR_l\\
\bigoplus\\
\bigoplus\limits_{0\leq l<N} V_{l+1}\otimes\CR_l\otimes\CO_\CS(-1,0)\\
\bigoplus\\
\bigoplus\limits_{0\leq l<N} W_l\otimes\CR_l
\end{array}
\stackrel{D}{\longrightarrow}
\bigoplus\limits_{0\leq l<N} V_{l+1}\otimes\CR_l
\end{equation}
Here the ``matrix coefficients'' of $C,D$ are as follows:
$_{ll}C^V_{ll}=tA_l-z;\qquad _{ll}C^V_{l,l-1}=-\xi_l;\qquad
_{NN}C^V_{10}=xB_0$, and $_{ll}C^V_{l+1,l}=B_l$ for $0<l<N$; furthermore,
$_{NN}C^W_{00}=txq_0$, and $_{ll}C^W_{ll}=tq_l$ for $0<l<N$. Furthermore,
$_{l+1,l}D^V_{l+1,l}=tA_l-z;\qquad _{ll}D^V_{l,l-1}=\xi_l;\qquad
_{NN}D^V_{10}=-xB_0$, and $_{ll}D^V_{l+1,l}=-B_l$ for $0<l<N$; furthermore,
$_{ll}D^W_{l+1,l}=p_{l+1}$.
We have used some evident shortcuts to simplify the notations, e.g.
$_{NN}C^V_{10}=xB_0:=B_0\otimes x\otimes1\in
\Hom(V_0,V_1)\otimes\Hom(\CR_N,\CR_0)\otimes\Hom_\bC(\CO_\bC(-1),\CO_\bC(-1))$.

\ssec{inverse}{Inverse construction}
Conversely, given a torsion free sheaf $\CF$ on $\CS_N$ framed at
$\CalD_\infty$, and $0\leq l<N$, we have (cf.~Section~5 of~\cite{BGK})
$H^0(\CS,\CR^*_l\otimes\CF(-1,0))=H^0(\CS,\CR^*_l\otimes\CF(0,-N))=
H^0(\CS,\CF(-1,-N))=0,\
H^2(\CS,\CR^*_l\otimes\CF(-1,0))=H^2(\CS,\CR^*_l\otimes\CF(0,-N))=
H^2(\CS,\CF(-1,-N))=0$. Furthermore, for $0<l<N,\
H^1(\CS,\CR^*_l\otimes\CF(-1,0))\simeq H^1(\CS,\CR^*_l\otimes\CF(0,-N))
\simeq V_l$, and $H^1(\CS,\CF(-1,0))\simeq H^1(\CS,\CF(0,-N))\simeq
H^1(\CS,\CF(-1,-N))\simeq V_0$.
Furthermore, for $0\leq l<N-1$, we have a canonical exact sequence
$$0\to H^0(\CS,\CR^*_l\otimes\CF)\to W_l\stackrel{p_{l+1}}{\longrightarrow}
V_{l+1}\to H^1(\CS,\CR^*_{l+1}\otimes\CF)\to0,$$ and also
$$0\to H^0(\CS,\CR^*_{N-1}\otimes\CF)\to W_{N-1}\stackrel{p_0}{\longrightarrow}
V_0\to H^1(\CS,\CF)\to0,$$
The Beilinson spectral sequence for $\CF$ takes the form 
$(E_1^{i,j})_{i=-2,-1,0}=$
$$\bigoplus\limits_{0<l\leq N}\Ext^j(\CR_l(1,0),\CF)\otimes\CR_l(-1,0)\to
\begin{array}{c}
\bigoplus\limits_{0<l\leq N}\Ext^j(\CR_l,\CF)\otimes\CR_l\\ \bigoplus\\
\bigoplus\limits_{0\leq l<N}\Ext^j(\CR_{l+1}(1,0),\CF)\otimes\CR_l(-1,0)
\end{array}
\to \bigoplus\limits_{0\leq l<N}\Ext^j(\CR_{l+1},\CF)\otimes\CR_l,$$
that is
$$\xymatrix{
\bigoplus\limits_{0<l\leq N}V_l\otimes\CR_l(-1,0) \ar[r]
\ar@{-->}[drr]_{d_2^{-2,1}}
& \bigoplus\limits_{0<l\leq N}V_l\otimes\CR_l\oplus
\bigoplus\limits_{0\leq l<N}V_{l+1}\otimes\CR_l(-1,0) \ar[r]
& \bigoplus\limits_{0\leq l<N}H^1(\CR^*_{l+1}\otimes\CF)\otimes\CR_l\\
&& \bigoplus\limits_{0\leq l<N}H^0(\CR^*_l\otimes\CF)\otimes\CR_l
}$$
Finally, we can replace $H^1(\CR^*_{l+1}\otimes\CF)=\on{Coker}p_{l+1}$
(resp. $H^0(\CR^*_l\otimes\CF)=\on{Ker}p_{l+1}$) by $V_{l+1}$ (resp. $W_l$),
and lift the differential $d_2^{-2,1}:\ E_2^{-2,1}\to E_2^{0,0}$ to a morphism
$V_l\otimes\CR_l(-1,0)\to W_l\otimes\CR_l$. Replacing the spectral sequence
with the total complex we obtain the ADHM description~\refe{mon} of $\CF$.

\ssec{PNN1}{Monad for the stack $\CS'_N$}
We also consider the following version of the above construction.
Let $\CS'=\CS'_N$ be the stacky weighted projective plane $\BP^2(N,N,1)$.
More precisely, we consider the affine 3-space $\BA^3$ with coordinates
$(z_0,z_1,z_2)$, and with the action of $\BC^*$ given by
$c(z_0,z_1,z_2)=(c^Nz_0,c^Nz_1,cz_2)$. We define $\CS':=(\BA^3\setminus0)/\BC^*$.
We define $\ell\subset\CS'$ as the hyperplane $z_2=0$ (all the
points of this line have automorphism group $\BZ/N\BZ$), and we define
$\ell_\infty\subset\CS'$ as the hyperplane $z_0=0$. Note that $\ell_\infty\simeq
\CX$. We denote $\CO_{\CS'}(l\ell)$
by $\CO(l)$ for short; note that $\CO_{\CS'}(\ell_\infty)\simeq\CO(N)$.
Let $\CP'_{\unl{d}}$ be the moduli space of torsion free sheaves of degree
$\unl{d}$ on $\CS'$ framed at $\ell_\infty$, i.e. such that $\CF|_{\ell_\infty}=
\CF_\infty:=W_0\otimes\CO_\CX\oplus W_1\otimes\CO_\CX(-1)\oplus\ldots\oplus
W_{N-1}\otimes\CO_\CX(1-N)$.

Since $\CS-\CalD_\infty\simeq\CS'-\ell_\infty$, and the framings at infinities
match, we have an identification $\CP_{\unl{d}}\simeq\CP'_{\unl{d}}$. We
describe the resulting isomorphism $\fM_{\unl{d}}\iso\CP'_{\unl{d}}$.
It sends a representative
$(A_\bullet,B_\bullet,p_\bullet,q_\bullet)$ to the middle cohomology of
the following monad of vector bundles on $\CS'$:
\eq{mon'}
\bigoplus_{0<l\leq N} V_l(-l)
\stackrel{C}{\longrightarrow}
\begin{array}{c}
\bigoplus\limits_{0<l\leq N} V_l(1-l)\\
\bigoplus\\
\bigoplus\limits_{0<l\leq N} V_l(N-l)\\
\bigoplus\\
\bigoplus\limits_{0\leq l<N} W_l(-l)
\end{array}
\stackrel{D}{\longrightarrow}
\bigoplus\limits_{0<l\leq N} V_l(N+1-l)
\end{equation}
Here the ``matrix coefficients'' of $C,D$ are as follows:
$-z_2:\ V_l(-l)\to V_l(1-l);\ z_0B_0:\ V_N(-N)\to V_1(0);\
B_l:\ V_l(-l)\to V_{l+1}(-l);\ z_1-z_0A_l:\ V_l(-l)\to V_l(N-l);\
z_0q_0:\ V_N(-N)\to W_0;\ q_l:\ V_l(-l)\to W_l(-l)$; furthermore,
$z_1-z_0A_l:\ V_l(1-l)\to V_l(N+1-l);\ z_2:\ V_l(N-l)\to V_l(N+1-l);\
-z_0B_0:\ V_N(0)\to V_1(N);\ -B_l:\ V_l(N-l)\to V_l(N-l)\to V_{l+1}(N-l);\
z_0p_{l+1}:\ W_l(-l)\to V_{l+1}(N-l)$.

\ssec{rot}{Rotation and the inverse construction}
Conversely, given a torsion free sheaf $\CF$ on $\CS'$ framed at $\ell_\infty$,
and $l=1,\ldots,N$, we have $H^0(\CS',\CF(l-N-1))=H^2(\CS',\CF(l-N-1))=0$,
and $V_l=H^1(\CS',\CF(l-N-1))$. The endomorphisms $A_l$ arise from the action
of $z_1\in\Gamma(\CS',\CO(N))$, and $B_l$ arises from the action of
$z_2\in\Gamma(\CS',\CO(1))$. More precisely, for $l\in\BZ$, we have
the morphism $z_2:\ \CF(l-N-1)\to\CF(l-N)$ which induces $B_l:\ V_l\to V_{l+1}$
for $1\leq l\leq N-1$, and also $z_2:\ H^1(\CS',\CF(-N-1))\to H^1(\CS',\CF(-N))
=V_1$. However, the short exact $0\to\CF(-N-1)\stackrel{z_0}{\longrightarrow}
\CF(-1)\to\CF_\infty(-1)\to0$ gives rise to the long exact sequence of
cohomology including $z_0:\ H^1(\CS',\CF(-N-1))\iso H^1(\CS',\CF(-1))=V_N$.
So we define $B_0:\ V_N\to V_1$ as the composition $z_2z_0^{-1}$.

Furthermore, we define $A_0:\ V_N\to V_N$ as the composition $z_1z_0^{-1}$.

Furthermore, the short exact sequence
$0\to\CF(-N)\stackrel{z_0}{\longrightarrow}\CF\to\CF_\infty\to0$
gives rise to the long exact sequence of cohomology including
$W_0=H^0(\CX,\CF_\infty)\to H^1(\CS',\CF(-N))=V_1$. We define $p_1$ as this
latter map $W_0\to V_1$.

Furthermore, for $0\leq l<N$, the short exact sequence
$$0\to\CF(l-2N-1)\stackrel{z_0}{\longrightarrow}\CF(l-N-1)\to\CF_\infty(l-N-1)
\to0$$ gives rise to the long exact sequence of cohomology including
$H^1(\CS',\CF(l-N-1))\to H^1(\CX,\CF_\infty(l-N-1))=W_l\oplus W_{l+1}\oplus
\ldots\oplus W_{N-1}$. For $0<l<N$, we define $q_l:\ V_l=H^1(\CS',\CF(l-N-1))\to
W_l$ as the direct summand of the above morphism. For $l=0$, we define
$q_0: V_N=H^1(\CS',\CF(-1))\to W_0$ as the composition of the direct summand
of the above morphism with $z_0^{-1}:\ V_N=H^1(\CS',\CF(-1))\to
H^1(\CS',\CF(-N-1))$.

It remains to define $A_l,\ l\ne0$, and $p_l,\ l\ne1$. To this end, we define
the {\em rotation} $\rho\unl{d}$ as follows: $\rho d_l:=d_{l+1},\ l\in\BZ/N\BZ$.
We have a natural {\em rotation} isomorphism
$R:\ \fM_{\unl{d}}\iso\fM_{\rho\unl{d}}$, taking the quiver data
$(V_\bullet,W_\bullet,A_\bullet,B_\bullet,p_\bullet,q_\bullet)$ to
$(V_{\bullet-1},W_{\bullet-1},A_{\bullet-1},B_{\bullet-1},p_{\bullet-1},
q_{\bullet-1})$. We define the corresponding isomorphism
$R:\ \CP'_{\unl{d}}\iso\CP'_{\rho\unl{d}}$ presently.

Given a framed torsion free sheaf $\CF$ on $\CS'$, we define
$\CG=R(\CF)$ as the kernel of the natural projection $\CF(1)\twoheadrightarrow
\imath_*W_0(1)$. Here $\imath$ stands for the closed embedding
$\CX\simeq\ell_\infty\hookrightarrow\CS'$. We have an exact sequence
$0\to W_0(1-N)\to\imath^*\CG\stackrel{r}{\longrightarrow}\imath^*\CF(1)\to
W_0(1)\to0,$
and the morphism $r$ factors as the composition
$\imath^*\CG\to W_1\oplus W_2(-1)\oplus\ldots\oplus W_{N-1}(2-N)\to
\imath^*\CF(1)$. Since for any $l=1,\ldots,N-1$ we have
$\Ext^1_\CX(W_l(1-l),W_0(1-N))=0$, we conclude that
$\CG|_{\ell_\infty}\simeq W_1\otimes\CO_\CX\oplus W_2(-1)\oplus\ldots\oplus
W_{N-1}(2-N)\oplus W_0(1-N)$.

Furthermore, the long exact cohomology sequence arising from the short exact
sequence $0\to\CG(l-N-1)\to\CF(l-N)\to\imath_*W_0(l-N)\to0$ implies
$H^1(\CS',\CG(l-N-1))=V_{l+1}$ for $0<l<N$. Also, the long exact cohomology
sequence arising from the short exact sequence
$0\to\CF(-N)\to\CG(-1)\to\imath_*(W_1(-1)\oplus\ldots\oplus W_{N-1}(1-N))\to0$
implies $H^1(\CS',\CG(-1))=H^1(\CS',\CF(-N))=V_1$. Finally, it is clear that
$R^N=\Id:\ \CP'_{\unl{d}}\to\CP'_{\unl{d}}$.

Returning to the definition of $A_l,p_l$, we set $A_l:=R^{-l}A_0R^l,\
p_l:=R^{1-l}p_1R^{l-1}$.

\ssec{gammak}{The action of $\Gamma_k$}
Let $\Gamma_k\simeq\BZ/k\BZ$ (resp. $\Gamma_{kN}\simeq\BZ/kN\BZ$)
be the group of $k$-th (resp. $kN$-th) roots of unity, with
generator $\zeta_k$ (resp. $\zeta_{kN}$). We have a surjection
$\Gamma_{kN}\twoheadrightarrow\Gamma_k,\ \zeta_{kN}\mapsto\zeta_k$.
Recall the stack $\CX=\CX_N=(\BA^2\setminus 0)/\BC^*$ where the action of
$\BC^*$ is given by $c(z_1,z_2)=(c^Nz_1,cz_2)$. The group $\Gamma_k$ acts
on $\CX_N$ as follows: $\zeta_k(z_1,z_2)=(z_1,\zeta_kz_2)$. The quotient
stack $\CX_N/\Gamma_k$ is denoted by $_k\CX_N$. We also have a 
$\Gamma_k$-equivariant morphism 
$\nu:\ \CX_N\to\CX_N,\ \nu(z_1,z_2)=(z_1^k,z_2^k)$. It factors through
$\CX_N\to\ _k\CX_N\stackrel{\Theta}{\longrightarrow}\CX_N$.

In the local coordinates of~\refss{S_N}, $\CX_N$ is glued from the affine
line $(\BA^1,y_1)$ with coordinate $y_1$, and $(\BA^1,y_2)/\Gamma_N$: both
$\BA^1\setminus\{y_1=0\}$ and $(\BA^1\setminus\{y_2=0\})/\Gamma_N$ coincide with
${\mathbb G}_m$, and we glue the charts with the help of $y_1=y_2^{-N}$.
Now $_k\CX_N$ is glued from $(\BA^1,y_1)/\Gamma_k$ and $(\BA^1,y_2)/\Gamma_{kN}$
with the help of $y_1^k=y_2^{-kN}$. Note that the group $\Gamma_k$ acts on the
chart $(\BA^1,y_2)/\Gamma_N$ as the quotient of $\Gamma_{kN}$ by the normal
subgroup $\Gamma_N\subset\Gamma_{kN}$.

The group $\Gamma_k$ acts on $\CS_N=\bC\times\CX_N$ via the second factor,
and we denote $\CS_N/\Gamma_k$ by $_k\CS_N$. By an abuse of notation we
denote by $\Theta$ the morphism $\on{id}\times\Theta:\ _k\CS_N\to\CS_N$.
The corresponding action of $\Gamma_k$ on $\bS$ is as follows: 
$\zeta_k(z:t,y:x)=(z:t,\zeta_ky:x)$. 
The corresponding action of $\Gamma_k$ on $\CS'$ is given by
$\zeta_k(z_0,z_1,z_2)=(z_0,z_1,\zeta_kz_2)$.

The group $\Gamma_k$ acts on the moduli space $\CP_{\unl{d}}$ of parabolic
sheaves trivialized at infinity via its action on $\bS$ and the {\em trivial}
action on the trivialization at infinity. The fixed point variety
$\CP_{\unl{d}}^{\Gamma_k}$ has various connected components, and we are going to
describe them in quiver terms. To this end note that $\Gamma_k$ acts on
the moduli space $\CP_{\unl{d}}=\fM_{\unl{d}}$
(resp. $\CP'_{\unl{d}}=\fM_{\unl{d}}$) of torsion free sheaves on $\CS$
(resp. $\CS'$) framed at $\CalD_\infty$ (resp. $\ell_\infty$) via its action on 
$\CS$ (resp. $\CS'$) and the {\em trivial} action on the framing. We have
$\CP_{\unl{d}}^{\Gamma_k}=\fM_{\unl{d}}^{\Gamma_k}=(\CP'_{\unl{d}})^{\Gamma_k}$. 
According to~\refss{local} (resp.~\refss{PNN1})
$\Gamma_k$ acts on $\fM_{\unl{d}}$ as follows:
$\zeta_k(A_\bullet,B_\bullet,p_\bullet,q_\bullet)=
(A_\bullet,\zeta_kB_\bullet,\zeta_kp_\bullet,q_\bullet)$. To formulate
the conclusion we consider the representations of the following quiver $Q^k$:

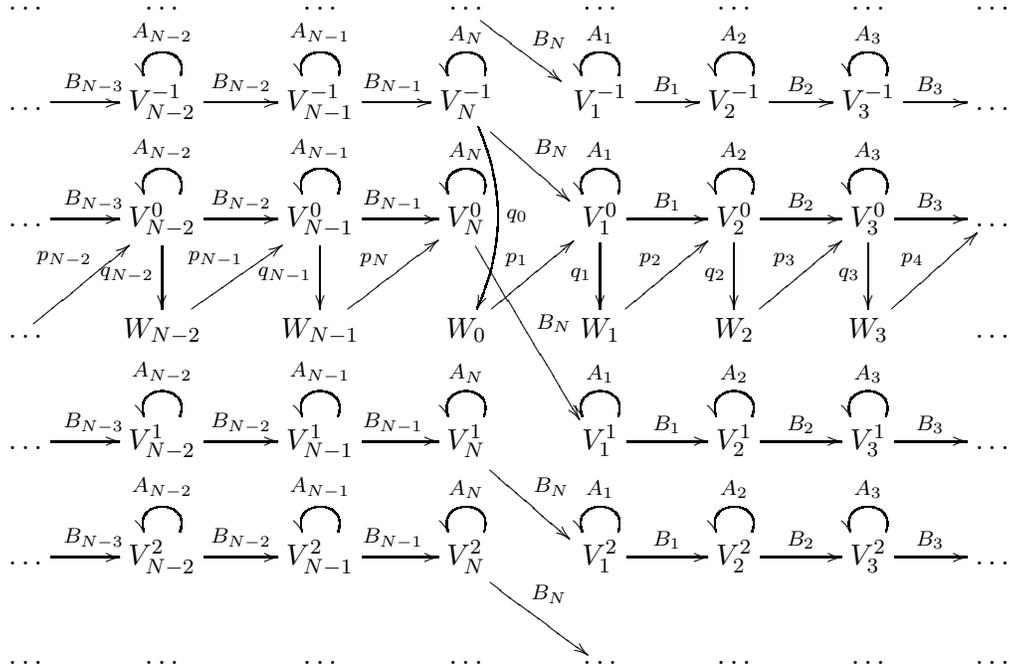
$$\xymatrix{
\ldots & \ldots & \ldots & \ldots \ar@<-.5ex>[dr]^{B_N} & \ldots
& \ldots & \ldots & \ldots\\
\ldots \ar[r]^{B_{N-3}}
& V_{N-2}^{-1} \ar@(ur,ul)[]_{A_{N-2}} \ar[r]^{B_{N-2}}
& V_{N-1}^{-1} \ar@(ur,ul)[]_{A_{N-1}} \ar[r]^{B_{N-1}}
& V_N^{-1} \ar@<-.5ex>[dr]^{B_N} \ar@(ur,ul)[]_{A_N} \ar@/^1pc/[dd]^{q_0}
& V_1^{-1} \ar[r]^{B_1} \ar@(ur,ul)[]_{A_1}
& V_2^{-1} \ar@(ur,ul)[]_{A_2} \ar[r]^{B_2}
& V_3^{-1} \ar@(ur,ul)[]_{A_3} \ar[r]^{B_3} &\ldots\\
\ldots \ar[r]^{B_{N-3}}
& V_{N-2}^0 \ar@(ur,ul)[]_{A_{N-2}} \ar[r]^{B_{N-2}} \ar[d]_{q_{N-2}}
& V_{N-1}^0 \ar@(ur,ul)[]_{A_{N-1}} \ar[r]^{B_{N-1}} \ar[d]_{q_{N-1}}
& V_N^0 \ar@<-.5ex>[ddr]^{B_N} \ar@(ur,ul)[]_{A_N}
& V_1^0 \ar[r]^{B_1} \ar[d]_{q_1} \ar@(ur,ul)[]_{A_1}
& V_2^0 \ar@(ur,ul)[]_{A_2} \ar[r]^{B_2} \ar[d]_{q_2}
& V_3^0 \ar@(ur,ul)[]_{A_3} \ar[r]^{B_3} \ar[d]_{q_3}
&\ldots\\
\ldots \ar[ur]^{p_{N-2}}
& W_{N-2} \ar[ur]^{p_{N-1}}
& W_{N-1} \ar[ur]^{p_N}
& W_0 \ar[ur]^{p_1} & W_1 \ar[ur]^{p_2} & W_2 \ar[ur]^{p_3}
& W_3 \ar[ur]^{p_4}
&\ldots\\
\ldots \ar[r]^{B_{N-3}}
& V_{N-2}^1 \ar@(ur,ul)[]_{A_{N-2}} \ar[r]^{B_{N-2}}
& V_{N-1}^1 \ar@(ur,ul)[]_{A_{N-1}} \ar[r]^{B_{N-1}}
& V_N^1 \ar@<-.5ex>[dr]^{B_N} \ar@(ur,ul)[]_{A_N}
& V_1^1 \ar[r]^{B_1} \ar@(ur,ul)[]_{A_1}
& V_2^1 \ar@(ur,ul)[]_{A_2} \ar[r]^{B_2}
& V_3^1 \ar@(ur,ul)[]_{A_3} \ar[r]^{B_3} &\ldots\\
\ldots \ar[r]^{B_{N-3}}
& V_{N-2}^2 \ar@(ur,ul)[]_{A_{N-2}} \ar[r]^{B_{N-2}}
& V_{N-1}^2 \ar@(ur,ul)[]_{A_{N-1}} \ar[r]^{B_{N-1}}
& V_N^2 \ar@<-.5ex>[dr]^{B_N} \ar@(ur,ul)[]_{A_N}
& V_1^2 \ar[r]^{B_1} \ar@(ur,ul)[]_{A_1}
& V_2^2 \ar@(ur,ul)[]_{A_2} \ar[r]^{B_2}
& V_3^2 \ar@(ur,ul)[]_{A_3} \ar[r]^{B_3} &\ldots\\
\ldots & \ldots & \ldots & \ldots & \ldots
& \ldots & \ldots & \ldots
}$$
Here the lower indices of $V$ run through $\BZ/N\BZ$, while the upper
indices run through $\BZ/k\BZ$. The relations are as follows:
$0=A_1B_N-B_NA_0+p_1q_0:\ V_N^{-1}\to V_1^0;\
0=A_{l+1}B_l-B_lA_l+p_{l+1}q_l:\ V_l^0\to V_{l+1}^0$ for $1\leq l\leq N-1$, and
$0=A_1B_N-B_NA_0:\ V_N^r\to V_1^{r+1}$ for $r\ne-1$, and
$0=A_{l+1}B_l-B_lA_l:\ V_l^r\to V_{l+1}^r$ in the remaining cases.

We set $d_l^r:=\dim(V_l^r)$, and we denote by $\widetilde{\unl{d}}$ the
collection of positive integers $(d_l^r)_{l\in\BZ/N\BZ}^{r\in\BZ/k\BZ}$.
We set $\unl{d}(\widetilde{\unl{d}}):=(d_1,\ldots,d_N)$ where
$d_l=\sum_{r\in\BZ/k\BZ}d_l^r$. We denote by $\sM_{\widetilde{\unl{d}}}$ the
scheme of representations of $Q^k$ of dimension $\widetilde{\unl{d}}$.
We denote by $G_{\widetilde{\unl{d}}}$ the group
$\prod_{l\in\BZ/N\BZ}^{r\in\BZ/k\BZ}\on{GL}(V_l^r)$; it acts naturally on
$\sM_{\widetilde{\unl{d}}}$. We denote by $\fZ_{\widetilde{\unl{d}}}$ the
categorical quotient $\sM_{\widetilde{\unl{d}}}/\!/G_{\widetilde{\unl{d}}}$.
Furthermore, we consider an open subscheme $\sM^s_{\widetilde{\unl{d}}}\subset
\sM_{\widetilde{\unl{d}}}$ formed by all the {\em stable} representations of
$Q^k$, i.e. those $(A_\bullet,B_\bullet,p_\bullet,q_\bullet)\in
\sM_{\widetilde{\unl{d}}}$ such that there is no proper graded subspace
$'V^\bullet_\bullet\subset V^\bullet_\bullet$ stable under $A_\bullet,B_\bullet$
and containing $p(W_\bullet)$. The action of $G_{\widetilde{\unl{d}}}$ on
$\sM^s_{\widetilde{\unl{d}}}$ is free, and we consider the GIT quotient
$\fM_{\widetilde{\unl{d}}}=\sM^s_{\widetilde{\unl{d}}}/G_{\widetilde{\unl{d}}}$.
Note that $\fM_{\widetilde{\unl{d}}}$ is nonempty iff
$d_N^0\geq d_1^1\geq d_2^1\geq\ldots\geq d_N^1\geq d_1^2\geq d_2^2\geq\ldots\geq
d_{N-2}^{-1}\geq d_{N-1}^{-1}\geq d_N^{-1}$.

The above considerations imply the following

\prop{invarianty}
The fixed point variety $\CP_{\unl{d}}^{\Gamma_k}$ is a union of connected
components isomorphic to $\fM_{\widetilde{\unl{d}}}$
(to be denoted by $\CP_{\widetilde{\unl{d}}}$), over all collections
$\widetilde{\unl{d}}$ such that $\unl{d}(\widetilde{\unl{d}})=\unl{d}$.
\eprop

\ssec{dir}{Direct image}
Given $\unl{d}=(d_1,\ldots,d_N)$ we consider $\widetilde{\unl{d}}=
\widetilde{\unl{d}}(\unl{d})$ such that
$d_l^0=d_l$ for any $1\leq l\leq N$, and $d_l^r=d_N$ for any $1\leq l\leq N$
and $r\ne0$ (note that $\unl{d}(\widetilde{\unl{d}}(\unl{d}))=
(d_1+(k-1)d_N,\ldots,d_N+(k-1)d_N)=:\unl{d}+(k-1)d_N$). 
Then it is easy to see that $\fM_{\widetilde{\unl{d}}}\simeq
\fM_{\unl{d}}$. In effect, all the maps $B_l$ except for
$B_N:\ V_N^{-1}\to V_1^0$, and the ones in the 0th row, have to be
isomorphisms intertwining the corresponding endomorphisms $A_l$ and $A_{l+1}$.

Geometrically, the isomorphism $\CP_{\widetilde{\unl{d}}}\iso
\CP_{\unl{d}}$ has the following explanation. We have an evident projection
$\psi:\ \bS\to\bS/\!/\Gamma_k\simeq\bS$ (the categorical quotient). A
$\Gamma_k$-fixed point of $\CP_{\unl{d}+(k-1)d_N}$ is represented by a
$\Gamma_k$-equivariant torsion free parabolic sheaf $\CF_\bullet$ on $\bS$.
Then $\psi_*\CF_\bullet$ carries a fiberwise action of $\Gamma_k$, and
$(\psi_*\CF_\bullet)^{\Gamma_k}$ is a torsion free parabolic sheaf on $\bS$,
trivialized at infinity. Its class in $\CP_{\unl{d}}$ is the image of $\CF_\bullet$
under the above isomorphism. For this reason, somewhat abusing notation, we
will denote this isomorphism by $\psi_*^{\Gamma_k}$.

Alternatively, thinking of $\Gamma_k$-equivariant torsion free parabolic 
sheaves on $\bS$ trivialized at infinity as of $\Gamma_k$-equivariant
torsion free sheaves on $\CS_N$ framed at infinity, that is torsion free
sheaves on $_k\CS_N$ framed at infinity (i.e. whose restriction to
$_k\CalD_\infty=\infty_\bC\times\ _k\CX_N\cup\bC\times\infty_{_k\CX_N}$ is equipped
with an isomorphism to $\CO_{_k\CX_N}\oplus\CO_{_k\CX_N}(-1\cdot0_{_k\CX_N})\oplus
\ldots\oplus\CO_{_k\CX_N}((1-N)\cdot0_{_k\CX_N})$ on $\infty_\bC\times\ _k\CX_N$,
and to $\CO_\bC^N$ on $\bC\times\infty_{_k\CX_N}$) we see that their 
isomorphism classes correspond bijectively to the isomorphism classes of 
stable representations of $Q^k$ (an argument entirely similar 
to~\refss{local} and~\refss{inverse}). Then $\psi_*^{\Gamma_k}$ is nothing
but $\Theta_*$ (notation of~\refss{gammak}).

For an arbitrary $\widetilde{\unl{d}}$, we consider $\unl{d}=(d_1,\ldots,d_N):=
(d_1^0,\ldots,d_{N-1}^0,d_N^{-1})$. Then we still have a morphism
$\Theta_*=\psi_*^{\Gamma_k}:\ \CP_{\widetilde{\unl{d}}}\to\CP_{\unl{d}}$
which is not necessarily an isomorphism. 
Going through the inverse constructions of~\refss{local} and~\refss{inverse} 
one arrives at the following description of
$\Theta_*=\psi_*^{\Gamma_k}:\ \fM_{\widetilde{\unl{d}}}\to\fM_{\unl{d}}$ in quiver terms.
We have $V'_N:=V_N^{-1},\ V'_l:=V_l^0$ for $1\leq l\leq N-1$.
Furthermore, we have $A'_l:=A_l$ for $1\leq l\leq N$, and
$B'_{N-1}$ is the composition of all $B$'s going from $V_{N-1}^0$ to $V_N^0$,
then to $V_1^1$, and all the way through to $V_N^{-1}$; while all the other
$B'_l$ coincide with the corresponding $B_l$.
Finally, $q'_l$ coincides with the corresponding $q_l$, and for
$1\leq l\leq N-1$ the map $p'_l$ coincides with the corresponding $p_l$; while
$p'_N$ is the composition of $p_N$ with all the $B$'s going from $V_N^0$ to
$V_1^1$, and all the way through to $V_N^{-1}$.

The morphism $\psi_*^{\Gamma_k}:\ \fM_{\widetilde{\unl{d}}}\to\fM_{\unl{d}}$ induces the
morphism $\psi^k:\ \fZ_{\widetilde{\unl{d}}}\to\fZ_{\unl{d}}$ from the affinization
$\fZ_{\widetilde{\unl{d}}}$ of $\fM_{\widetilde{\unl{d}}}$ to the affinization
$\fZ_{\unl{d}}$ of $\fM_{\unl{d}}$.

\ssec{def}{Defect}
Let us give a geometric explanation of what is so special about the
components $\CP_{\widetilde{\unl{d}}},\ \widetilde{\unl{d}}=
\widetilde{\unl{d}}(\unl{d})$ considered in~\refss{dir}. Namely, they are
the only components of the fixed point variety $\CP^{\Gamma_k}_{\unl{d}}$ which
contain the {\em nonempty} open subset formed by the $\Gamma_k$-equivariant
{\em locally free} parabolic sheaves. Here {\em locally free} means
locally free after forgetting the $\Gamma_k$-equivariant structure.

For an arbitrary $\Gamma_k$-equivariant torsion free parabolic sheaf
$\CF_\bullet$, there is a notion
of the {\em saturation} $\widehat\CF_\bullet$ (a locally free parabolic sheaf
containing $\CF_\bullet$, such that the quotient has a zero-dimensional
support). The global sections of this quotient is a $\BZ/N\BZ$-graded
$\Gamma_k$-module $\on{def}(\CF_\bullet)$, the {\em defect} of $\CF_\bullet$.
The class $[\on{def}(\CF_\bullet)]$ of $\on{def}(\CF_\bullet)$
in the $K$-group of $\BZ/N\BZ$-graded
$\Gamma_k$-modules is represented by a collection $\widetilde{\unl{d}}$ of
integers. The class $[\on{def}(\CF_\bullet)]$ may vary throughout a connected
component of the fixed point variety $\CP^{\Gamma_k}_{\unl{d}}$. However,
its class $[\overline{\on{def}}(\CF_\bullet)]$ modulo the subgroup spanned
by all the collections of the sort $\widetilde{\unl{d}}(\unl{d}'),\
\unl{d}'\in\BZ^{\BZ/N\BZ}$, is constant throughout a connected component
$\CP_{\widetilde{\unl{d}}}$. Quite evidently, the class
$[\overline{\on{def}}(\CF_\bullet)]$ for
$\CF_\bullet\in\CP_{\widetilde{\unl{d}}}$, equals the class of
$\widetilde{\unl{d}}$. In particular, in order to have a locally free
parabolic sheaf $\CF_\bullet$ (i.e. the one with zero defect) in a component
$\CP_{\widetilde{\unl{d}}}$ it is necessary and sufficient that
$\widetilde{\unl{d}}$ be of the form $\widetilde{\unl{d}}(\unl{d})$ for some
$\unl{d}$.

\sec{zastava blowup}{Zastava for blown up plane}
The results of this section are strongly influenced by~\cite{NY}.

\ssec{den chain}{Dented chainsaw}
We consider the representations of the following {\em dented chainsaw quiver} $\hQ$

$$\xymatrix{
\ldots \ar[r]^{B_{N-3}}
& V_{N-2} \ar@(ur,ul)[]_{A_{N-2}} \ar[r]^{B_{N-2}} \ar[d]_{q_{N-2}}
& V_{N-1} \ar@(ur,ul)[]_{A_{N-1}} \ar[r]^{B_{N-1}} \ar[d]_{q_{N-1}}
& V_N \ar@<-.5ex>[r]_e
& V_0 \ar@<-.5ex>[l]_\delta \ar[r]^{B_0} \ar[d]_{q_0}
& V_1 \ar@(ur,ul)[]_{A_1} \ar[r]^{B_1} \ar[d]_{q_1}
& V_2 \ar@(ur,ul)[]_{A_2} \ar[r]^{B_2} \ar[d]_{q_2} &\ldots\\
\ldots \ar[ur]^{p_{N-2}} & W_{N-2} \ar[ur]^{p_{N-1}} & W_{N-1} \ar[ur]^{p_N} &
& W_0 \ar[ur]^{p_1} & W_1 \ar[ur]^{p_2} & W_2 \ar[ur]^{p_3} &\ldots
}$$
with relations $A_{l+1}B_l-B_lA_l+p_{l+1}q_l=0$ for any $1\leq l\leq N-2;\
A_1B_0-B_0e\delta+p_1q_0=0;\ \delta eB_{N-1}-B_{N-1}A_{N-1}+p_Nq_{N-1}=0$.
Here $\dim W_l=1,\
d_N:=\dim V_N=d_0:=\dim V_0,\ \dim V_l=d_l,\ l=1,\ldots,N-1$.
We denote by $\hsM_{\unl{d}}$ the scheme of representations of $\hQ$.
We denote by $\hG_{\unl{d}}$ the group $\prod_{0\leq l\leq N}\GL(V_l)$;
it acts naturally on $\hsM_{\unl{d}}$.
Performing the celebrated Crawley-Boevey trick, we identify all the lines
$W_l$ with, say $W_\infty$, so that $W_\infty$ is the source of all $p_l$,
and the target of all $q_l$. We will denote a typical representation of
$\hQ$ by $Y$.

\ssec{stab}{Stability conditions}
Following~\cite[Section 4(ii)]{NaM}
we consider the enhanced dimension vectors $\hd:=(d_0,d_1,d_2,
\ldots,d_{N-1},d_N)$, and $\td:=(d_0,d_1,d_2,
\ldots,d_{N-1},d_N,1)$ with one extra coordinate equal to $\dim W_\infty=1$.
We consider a vector $\zeta^\bullet=(\zeta_0,\zeta_1,\ldots,\zeta_{N-1},\zeta_N)$
where $\zeta_N=-1,\ \zeta_0=1,\ \zeta_l=0$ for $l=1,\ldots,N-1$.
Also, for $0<\varepsilon\ll1$ we consider $\zeta^-:=\zeta^\bullet-(\varepsilon,
\ldots,\varepsilon)$. We set $\zeta^-_\infty:=-\langle\zeta^-,\hd\rangle$,
and $\zeta^\bullet_\infty:=-\langle\zeta^\bullet,\hd\rangle$ where
$\langle\cdot,\cdot\rangle$ stands for the sum of products of coordinates
(the standard scalar product). Finally, we set
$\tzeta^-:=(\zeta^-,\zeta^-_\infty)$, and
$\tzeta^\bullet:=(\zeta^\bullet,\zeta^\bullet_\infty)$.

For a nonzero $\hQ$-submodule $Y'\subset Y$ of enhanced dimension $\td'$
(where the last coordinate may be either 1 or 0) we define the slope by
$$\theta^-(Y'):=\frac{\langle\tzeta^-,\td'\rangle}
{\langle(1,\ldots,1),\td'\rangle},\
\theta^\bullet(Y'):=\frac{\langle\tzeta^\bullet,\td'\rangle}
{\langle(1,\ldots,1),\td'\rangle}.$$
We say that a $\hQ$-module $Y$ is $\zeta^-$-semistable (resp.
$\zeta^\bullet$-semistable) if for any nonzero submodule $Y'\subset Y$
we have $\theta^-(Y')\leq\theta^-(Y)$ (resp.
$\theta^\bullet(Y')\leq\theta^\bullet(Y)$). We say $Y$ is
$\zeta^-$-stable (resp. $\zeta^\bullet$-stable) if the inequality is strict
unless $Y'=Y$. Note that $\zeta^-$-stability is equivalent to
$\zeta^-$-semistability.

We define a scheme $\hfM_{\unl{d}}$ as the moduli space
of $\zeta^-$-semistable (equivalently, $\zeta^-$-stable) $\hQ$-modules.
By GIT, $\hfM_{\unl{d}}$ is the projective spectrum of the ring of
$\hG_{\unl{d}}$-semiinvariants in $\BC[\hsM_{\unl{d}}]$.
Furthermore, we define a scheme $\hfZ_{\unl{d}}$ as the moduli space of
$S$-equivalence classes of $\zeta^\bullet$-semistable $\hQ$-modules.
Since the stability condition $\zeta^\bullet$ lies on a wall of the chamber
containing $\zeta^-$, we have a projective morphism $\pi_{\zeta^\bullet,\zeta^-}:\
\hfM_{\unl{d}}\to\hfZ_{\unl{d}}$.

\ssec{parab}{Parabolic sheaves on blow-up}
We stick to the notations of~\cite{NY}. Namely, $\BP^2$ is the projective
plane with homogeneous coordinates $[z_0:z_1:z_2]$, and $\ell_\infty\subset\BP^2$
is the line ``at infinity'' given by the equation $z_0=0$. Furthermore,
$\hBP^2$ is the blow-up of $\BP^2$ ``at the origin'' (given by equations
$z_1=z_2=0$). It is the closed subvariety of $\BP^2\times\BP^1$ defined by
$\hBP^2=\{([z_0:z_1:z_2],[z:w]):\ z_1w=z_2z\}$. We denote by $E$ the exceptional
divisor in $\hBP^2$; we denote by $\ell_0\subset\hBP^2$ the proper transform
of the line $z_2=0$ in $\BP^2$; finally, by an abuse of notation,
we denote by $\ell_\infty\subset\hBP^2$
the proper transform of the line $\ell_\infty\subset\BP^2$.

We set $W:=W_1\oplus W_2\oplus\ldots\oplus W_{N-1}\oplus W_0$.
Given an $N$-tuple of nonnegative integers $\unl{d}=(d_0,\ldots,d_{N-1})$
we say that a {\em parabolic sheaf} $\CF_\bullet$ of degree $\unl{d}$ is an
infinite flag of torsion free coherent sheaves of rank $N$ on $\hBP^2:\
\ldots\subset\CF_{-1}\subset\CF_0\subset\CF_1\subset\ldots$ such that

(a) $\CF_{k+N}=\CF_k(\ell_0)$ for any $k\in\BZ$;

(b) $ch_1(\CF_k)=k[\ell_0]$ for any $k\in\BZ$: the first Chern classes are
proportional to the fundamental class of $\ell_0$;

(c) $ch_2(\CF_k)=d_i$ for $i\equiv k\pmod{N}$;

(d) $\CF_0$ is locally free at $\ell_\infty$ and trivialized at $\ell_\infty:\
\CF_0|_{\ell_\infty}=W\otimes\CO_{\ell_\infty}$;

(e) For $-N\leq k\leq0$ the sheaf $\CF_k$ is locally free at $\ell_\infty$,
and the quotient sheaves $\CF_k/\CF_{-N},\ \CF_0/\CF_k$ (both supported at
$\ell_0\subset\hBP^2$) are locally free at the point $\ell_0\cap\ell_\infty$;
moreover, the local sections of $\CF_k|_{\ell_\infty}$ are those sections of
$\CF_0|_{\ell_\infty}=W\otimes\CO_{\ell_\infty}$ which take value in
$W_1\oplus\ldots\oplus W_{k+N}\subset W$ at $\ell_0\cap\ell_\infty$.

One can show that the fine moduli space $\hCP_{\unl{d}}$ of degree $\unl{d}$
parabolic sheaves exists, and is a smooth connected quasiprojective variety
of dimension $2d_0+\ldots+2d_{N-1}$.

\th{queer blowup} {\em (A.~Kuznetsov)}
There is an isomorphism $\Xi:\ \hfM_{\unl{d}}\iso\hCP_{\unl{d}}$.
\eth

The proof occupies~\refss{hat S} --- \refss{ass}.

\ssec{hat S}{Stack $\hCS_N$}
We denote by $\ell'\subset\hBP^2$ the proper transform of the line
$z_1=0$ in $\BP^2$. We consider the open subvarieties
$U':=\hBP^2-\ell_0\simeq\ell'\times\BA^1$ with
coordinate $z=z_1z_2^{-1}$ along $\BA^1$,
and $U_0:=\hBP^2-\ell'\simeq\ell_0\times\BA^1$ with coordinate $w=z_2z_1^{-1}$
along $\BA^1$. Note that $\ell_0$ (resp. $\ell'$) in 
$\hBP^2\subset\BP^2\times\BP^1$ is cut out by the equation $w=0$ 
(resp. $z=0$).
We consider the ramified Galois covering
$\theta:\ \BA^1\to\BA^1,\ w=s^N$ with Galois group $\Gamma_N$.
We denote by $\theta:\ \widetilde{U}_0\to U_0$ the base change of this
covering under $U_0\to\BA^1$. The action of $\Gamma_N$ on
$\theta^{-1}(U_0\cap U')$ is free, and $\theta^{-1}(U_0\cap U')/\Gamma_N=
U_0\cap U'$. We define a 2-dimensional Deligne-Mumford stack $\hCS_N$ as
the result of gluing $U'$ and $\widetilde{U}_0/\Gamma_N$ over the common
open $U_0\cap U'$. Note that $\hCS_N$ is equipped with a projection
$\vartheta$ to $\hBP^2$ which is an isomorphism off $\ell_0$. The line in
$\hCS_N$ lying over $\ell_0$ will be denoted by $\ell\subset\hCS_N$; its
automorphism group is $\Gamma_N$. Since $N$ is fixed throughout the Section,
we will often omit the lower index $N$ to simplify the notations.

We also have a smooth morphism $\pi$ (a $\BP^1$-bundle) from $\hCS_N$ to the
1-dimensional stack $\CX_N$ of~\refss{S_N} such that
$\pi^{-1}(\infty_\CX)=\ell'$, and $\pi^{-1}(0_\CX)=\ell$.
A section of $\pi$ sending $\CX$ to $\ell_\infty\subset\hCS$ will be
denoted by $\imath$. We choose a section
$y_2$ of $\CO_\CX(1)$ with a simple zero at $0_\CX$, and a section $y_1$
of $\CO_\CX(N)$ with a simple zero at $\infty_\CX$ (in notations of~\refss{adhm}
and~\refss{S_N} we have $y_2^N=y,\ y_1=x$). We keep the same names
for the corresponding sections of $\CO_\hCS(\ell)$ and $\CO_\hCS(N\ell)$
constant along the fibers of $\pi$. Finally, we choose a section $x_2$
of $\CO_\hCS(\ell_\infty-N\ell)=\CO_\hCS(\ell_\infty-\ell')$
with a simple zero at $E$.

Here is an alternative toric description of $\hCS$. We consider
$\BA^4$ with coordinates $x_1,x_2,y_1,y_2$, and with an open subset 
$\widetilde{U}\subset\BA^4$ obtained by removing two planes:
$L_1=\{x_1,x_2,0,0\}$ and $L_2=\{0,0,y_1,y_2\}$. The torus $T_2=\BC^*\times\BC^*$
acts on $\widetilde U$ as follows: $(c_1,c_2)\cdot(x_1,x_2,y_1,y_2)=
(c_1c_2^Nx_1,c_1x_2,c_2^Ny_1,c_2y_2)$. We have $\hCS=\widetilde{U}/T_2$.
Note that $x_1\in\Gamma(\hCS,\CO(\ell_\infty))$ is an equation of $\ell_\infty$;
$x_2\in\Gamma(\hCS,\CO(E))$ is an equation of $E$;
$y_1\in\Gamma(\hCS,\pi^*\CO_\CX(N))$ is an equation of $\pi^{-1}(\infty_\CX)$;
$y_2\in\Gamma(\hCS,\pi^*\CO_\CX(1))$ is an equation of $\ell$.

According to~\cite{Bi}, there is a one-to-one correspondence between the
(torsion free, framed at $\ell_\infty$) sheaves on $\hCS$, and the
(torsion free, framed at $\ell_\infty$) sheaves on $\hBP^2$ with parabolic
structure along $\ell_0$. Thus $\hCP_{\unl{d}}$ is the moduli space of torsion
free sheaves of degree $\unl{d}$ on $\hCS$ framed at $\ell_\infty$.
More precisely, the framing at $\ell_\infty$ is an isomorphism
$\imath^*\CF\simeq\CF_\infty:=
W_0\otimes\CO_\CX\oplus W_1\otimes\CO_\CX(-1)\oplus
\ldots\oplus W_{N-1}\otimes\CO_\CX(-N+1)$.
%W_0\otimes\CO_\CX(-1)\oplus W_1\otimes\CO_\CX(-2)\oplus
%\ldots\oplus W_{N-1}\otimes\CO_\CX(-N)$.
For technical reasons, it will be more convenient for us to view
$\hCP_{\unl{d}}$ as the moduli space of twisted sheaves $\CG:=\CF(-\ell)$
with framing at $\ell_\infty:\ \imath^*\CG\simeq\CG_\infty:=
W_0\otimes\CO_\CX(-1)\oplus W_1\otimes\CO_\CX(-2)\oplus
\ldots\oplus W_{N-1}\otimes\CO_\CX(-N)$.

\ssec{plan}{Plan of the proof}
An exceptional collection of line bundles $\{\CO_\CX,\CO_\CX(1),\ldots,\CO_\CX(N)\}$ on
$\CX$ gives rise to an equivalence of the derived category $D(\CX)$ of coherent sheaves
on $\CX$ and the derived category $D(K)$ of representations of the following
quiver $K$:
\eq{xquiver}
\xymatrix{
V_0 \ar[r]^-{B_0} \ar@/_1pc/[rr]_{\delta} & \cdots \ar[r]^-{B_{N-1}} & V_N \\
}
\end{equation}

Since $\hCS$ is a $\BP^1$-bundle over $\CX$, it possesses the following exceptional
collection of line bundles: $\{\pi^*\CO_\CX(-N-1),\ldots,\pi^*\CO_\CX(-1),
(\pi^*\CO_\CX)(-\ell_\infty),\ldots,(\pi^*\CO_\CX(N))(-\ell_\infty)\}$. The derived category
$D(\hCS)$ is equivalent to the category of diagrams 
$\CG'(-N)\stackrel{\mu_0}{\longleftarrow}\CG''\stackrel{\mu_\infty}{\longrightarrow}\CG'$
where $\CG',\CG''\in D(\CX)\simeq D(K)$. For $\CG\in D(\hCS)$ we set $\CG'=\pi_*\CG,\
\CG''=\pi_*\CG(-\ell_\infty)$. We have an exact triangle 
\eq{32}
\ldots\to\pi^*\CG''(-\ell_\infty+N\ell)\stackrel{\mu}{\longrightarrow}\pi^*\CG'\to\CG\to
\pi^*\CG''(-\ell_\infty+N\ell)[1]\to\ldots
\end{equation}
By the adjointness and the projection formula, the morphism $\mu$ is the same as
the morphism $\mu_\infty\oplus\mu_0$ from $\CG''$ to 
$\CG'\otimes\pi_*\CO_\hCS(\ell_\infty-N\ell)=\CG'\oplus\CG'(-N)$.
The framing of $\CG$ at $\ell_\infty$ implies the existence of the following exact
triangle:
$$\ldots\to\CG''\stackrel{\mu_\infty}{\longrightarrow}\CG'\to\CG_\infty
\stackrel{\beta}{\longrightarrow}\CG''[1]\to\ldots$$
Now the dual exceptional collection of 
$\{\CO_\CX,\CO_\CX(1),\ldots,\CO_\CX(N)\}\in D(\CX)$ is 
$\{\CO_\CX(-1-N), \on{Coker}(\CO_\CX(-2-N)\hookrightarrow\CO_\CX(-1-N)),\ldots,
\on{Coker}(\CO_\CX(-2N)\hookrightarrow\CO_\CX(1-2N)), \CO_\CX(-N)\}$. Decomposing $\CG''$
with respect to this exceptional collection we obtain a representation of the 
quiver~\refe{xquiver}. The additional data $p_\bullet$ of the dented chainsaw quiver
correspond to the morphism $\beta$ (equivalently, $\mu_\infty$) above, 
while the additional data 
$(e,A_\bullet,q_\bullet)$ correspond to the morphism $\mu_0$ above. 
The property of 
$V_\bullet$ being vector spaces (as opposed to complexes of vector spaces) is 
equivalent to the property of $\CG$ being a perverse coherent sheaf with
torsion supported at the exceptional divisor $E$ and in codimension 2
(i.e. at finitely many points off $E$). Now $\zeta^\bullet$-semistability
is equivalent to the vanishing of torsion of $\CG$ at the generic point of $E$,
while the $\zeta^-$-semistability is equivalent to the latter vanishing plus
vanishing of the first cohomology of $\CG$, that is $\CG$ being a torsion
free sheaf.

\ssec{monad hat S}{Monad for the stack $\hCS_N$}
%Note that $\vartheta_*$ establishes an isomorphism
%$\Gamma(\hCS,\CO_\hCS(\ell_\infty-E))\simeq
%\Gamma(\hBP^2,\CO_{\hBP^2}(\ell_\infty-E)=\BC\langle z,w\rangle$, and
%$\Gamma(\hCS,\CO_\hCS(\ell_\infty))\simeq
%\Gamma(\hBP^2,\CO_{\hBP^2}(\ell_\infty)=\BC\langle z_0,z_1,z_2\rangle$.
Given a $\hQ$-module $Y$, we construct $\CG$ as 
the following monad of vector bundles on $\hCS$ (in cohomological degrees
$-1,0,1$):
\eq{mona}
\begin{array}{c}
%V_0(-\ell_\infty)\\
%\bigoplus\\
\bigoplus\limits_{l=0}^NV_l(-\ell_\infty-(l+1)\ell)
\end{array}
\stackrel{C}{\longrightarrow}
\begin{array}{c}
\bigoplus\limits_{l=1}^N V_l(-\ell_\infty-l\ell)\\
\bigoplus\\
V_N(-\ell_\infty-\ell)\\
\bigoplus\\
\bigoplus\limits_{l=0}^{N-1} W_l((-l-1)\ell)\\
\bigoplus\\
V_0((-N-1)\ell)\\
\bigoplus\\
\bigoplus\limits_{l=1}^{N-1}V_l((-l-1)\ell)
\end{array}
\stackrel{D}{\longrightarrow}
\bigoplus\limits_{l=1}^N V_l(-l\ell)
\end{equation}
%Here the middle term $\CM$ is an extension
%\eq{exten}
%0\to\bigoplus\limits_{l=0}^{N-1} W_l((-l-1)\ell)\oplus V_0((-N-1)\ell)\oplus
%\bigoplus\limits_{l=1}^{N-1}V_l((-l-1)\ell)
%\stackrel{\kappa}{\longrightarrow}
%\to\CM\to
%\stackrel{\varkappa}{\longrightarrow}
%\bigoplus\limits_{l=0}^N V_l(-\ell_\infty+(N-l)\ell)\to0
%\end{equation}
%given by an element
%$$\kappa\in\Ext^1_\hCS\left(\bigoplus\limits_{l=0}^N 
%V_l(-\ell_\infty+(N-l)\ell),
%\bigoplus\limits_{l=0}^{N-1} W_l((-l-1)\ell)\oplus V_0((-N-1)\ell)\oplus
%\bigoplus\limits_{l=1}^{N-1}V_l((-l-1)\ell)\right).$$
The morphisms $C,D$ are described as follows.

We introduce the complex $\CG''$ in the derived coherent category of $\CX$:
\eq{comple}
\bigoplus_{l=0}^NV_l(-N-l-1)\stackrel{\gamma''}{\longrightarrow}
\bigoplus_{l=1}^NV_l(-N-l)\oplus V_N(-N-1)
\end{equation}
Here the ``matrix coefficients'' of $\gamma''$ are as follows:
$y_2:\ V_l(-N-l-1)\to V_l(-N-l);\ B_l:\ V_l(-N-l-1)\to V_{l+1}(-N-l-1);\
-y_1:\ V_N(-2N-1)\to V_N(-N-1);\ \delta:\ V_0(-N-1)\to V_N(-N-1)$.

Alternatively, $\CG''$ is canonically quasiisomorphic to another complex
\eq{compl}
V_0(-N-1)\oplus\bigoplus_{l=1}^{N-1}V_l(-l-1)\stackrel{\gamma'}{\longrightarrow}
\bigoplus_{l=1}^NV_l(-l)
\end{equation}
where the ``matrix coefficients'' of $\gamma'$ are as follows:
$y_2:\ V_l(-l-1)\to V_l(-l);\ B_l:\ V_l(-l-1)\to V_{l+1}(-l-1);\
y_1B_0:\ V_0(-N-1)\to V_1(-1);\ y_2\delta:\ V_0(-N-1)\to V_N(-N)$.

The quasiisomorphism is given by
\eq{popytka}
\xymatrix{
\bigoplus_{l=0}^NV_l(-N-l-1) \ar[r]^-{\gamma} \ar[d]^\upsilon
& \bigoplus_{l=1}^NV_l(-N-l)\oplus V_N(-N-1) %\ar[dl]^{\varkappa} 
\ar[d]^{\upsilon'}\\
V_0(-N-1)\oplus\bigoplus_{l=1}^{N-1}V_l(-l-1) \ar[r]^-{\gamma'}
& \bigoplus_{l=1}^NV_l(-l)\\
}
\end{equation}
where the nonzero components of $\upsilon,\upsilon'$ are as follows:
$\on{id}:\ V_0(-N-1)\to V_0(-N-1);\ -y_1:\ V_l(-N-l-1)\to V_l(-l-1)$ 
(resp. $y_1:\ V_l(-N-l)\to V_l(-l)$) for $1\leq l\leq N-1;\
y_1:\ V_N(-2N)\to V_N(-N);\ y_2:\ V_N(-N-1)\to V_N(-N)$.

Finally, we are able to describe the morphisms $C,D$ of~\refe{mona}.
We have $C=C_1+C_2+C_3$ where $C_1$ is 
$$\pi^*\gamma''(-\ell_\infty+N\ell):\
\bigoplus\limits_{l=0}^NV_l(-\ell_\infty-(l+1)\ell)\to
\bigoplus\limits_{l=1}^N V_l(-\ell_\infty-l\ell)\oplus
V_N(-\ell_\infty-\ell),$$
%V_0(-\ell_\infty)\oplus\bigoplus\limits_{l=0}^{N-1}
%V_l(-\ell_\infty+(N-l-1)\ell)
%\to\bigoplus\limits_{l=0}^N V_l(-\ell_\infty+(N-l)\ell),$$ 
the matrix elements of $C_2:\ 
\bigoplus\limits_{l=0}^NV_l(-\ell_\infty-(l+1)\ell)\to
\bigoplus\limits_{l=0}^{N-1} W_l((-l-1)\ell)$ are
$x_1q_l\in\Hom(V_l(-\ell_\infty-(l+1)\ell),W_l((-l-1)\ell))$, i.e.
they correspond by the adjointness and projection formula to
$q_l\in\Hom_\CX(V_l(-l-1),W_l(-l-1))$.
Now $$C_3:\ \bigoplus\limits_{l=0}^NV_l(-\ell_\infty-(l+1)\ell)\to
V_0((-N-1)\ell)\oplus\bigoplus\limits_{l=1}^{N-1}V_l((-l-1)\ell)$$ 
corresponds by adjointness and projection formula to $\upsilon$ 
of~\refe{popytka} plus $$(0,A_1(-2),\ldots,A_{N-1}(-N),e(-N-1)):\
\bigoplus\limits_{l=0}^NV_l(-l-1)\to
V_0(-N-1)\oplus\bigoplus\limits_{l=1}^{N-1}V_l(-l-1).$$
In other words, $C_3=x_2\upsilon+x_1e+x_1\sum_{l=1}^{N-1}A_l$.

We have $D=D_1+D_2+D_3$ where 
$D_1:\ \bigoplus\limits_{l=1}^N V_l(-\ell_\infty-l\ell)\oplus
V_N(-\ell_\infty-\ell)\to\bigoplus\limits_{l=1}^N V_l(-l\ell)$ 
corresponds by adjointness and projection
formula to $\upsilon'$ of~\refe{popytka} plus
$$(-A_1(-1),\ldots,-A_{N-1}(-N+1),-\delta e(-N),B_0e(-1)):\ 
\bigoplus\limits_{l=1}^N V_l(-l)\oplus
V_N(-1)\to\bigoplus\limits_{l=1}^N V_l(-l).$$
In other words, $D_1=x_2\upsilon'-x_1\sum_{l=1}^{N-1}A_l-x_1\delta e+x_1B_0e$.
Now $D_2:\ \bigoplus\limits_{l=0}^{N-1} W_l((-l-1)\ell)\to
\bigoplus\limits_{l=1}^N V_l(-l\ell)$ is 
$-\pi^*\beta:\ \bigoplus\limits_{l=0}^{N-1} W_l((-l-1)\ell)
\to\bigoplus\limits_{l=1}^N V_l(-l\ell).$
The ``matrix elements'' of $\beta\in\Hom(\CG_\infty,\bigoplus_{l=1}^NV_l(-l))$
are $p_{l+1}:\ W_l(-l-1)\to V_{l+1}(-l-1)$. 
%(modulo the image of $\gamma'\circ?$).
Finally, $D_3:\ V_0((-N-1)\ell)\oplus\bigoplus\limits_{l=1}^{N-1}V_l((-l-1)\ell) 
\to\bigoplus\limits_{l=1}^N V_l(-l\ell)$ is $\pi^*\gamma'$.

%For the reader's convenience we collect some of the above morphisms in
%the following figure:
%\eq{xym}
%\xymatrix{
%& V_0(0)\oplus\bigoplus_{l=0}^{N-1}V_l(N-l-1) \ar[r]^-{\gamma}
%& \bigoplus_{l=0}^NV_l(N-l) \ar[dl]^{\kappa''_1}
%\ar@/_1pc/@{.>}[ddll]_{\kappa''_2}\\
%& \bigoplus_{l=0}^{N-1}W_l(-l-1) \ar[d]^\beta &\\
%V_0(-N-1)\oplus\bigoplus_{l=1}^{N-1}V_l(-l-1) \ar[r]^-{\gamma'}
%& \bigoplus_{l=1}^NV_l(-l) &\\
%}
%\end{equation}
%Here $\kappa''_1\circ\gamma=0$, while
%$\beta\in\Hom(\CG_\infty,\bigoplus_{l=1}^NV_l(-l))$ is only defined modulo the
%image of $\gamma'\circ?:\
%\Hom\left(\CG_\infty,V_0(-N-1)\oplus\bigoplus_{l=1}^{N-1}V_l(-l-1)\right)
%\to \Hom(\CG_\infty,\bigoplus_{l=1}^NV_l(-l))$.
%Finally, $\beta\circ\kappa''_1=\gamma'\circ\kappa''_2$, and
%$\kappa''_2\circ\gamma=0$.

\ssec{ic}{Inverse construction}
Conversely, given a torsion free sheaf $\CG$ on $\hCS$ with a framing
$\imath^*\CG\simeq\CG_\infty$ we have $H^0(\hCS,\CG(-\ell_\infty+l\ell))=
H^2(\hCS,\CG(-\ell_\infty+l\ell))=0$ for $0\leq l\leq N$, and we set
$V_l:=H^1(\hCS,\CG(-\ell_\infty+l\ell))$. Furthermore, we set
$\CG'':=\pi_*(\CG(-\ell_\infty))$, and $\CG':=\pi_*\CG$, so that
$V_l=H^1(\CX,\CG''(l))$.

Let $\Delta:\ \hCS\hookrightarrow\hCS\times\hCS$ stand for the diagonal
embedding.
We have the following exact triangles in the derived category of coherent
sheaves on $\hCS$:
\eq{conv1}
\ldots\to\CO_{\hCS\times_\CX\hCS}(-\ell_\infty,-\ell_\infty+N\ell)
\stackrel{\xi}{\longrightarrow}\CO_{\hCS\times_\CX\hCS}\to\Delta_*\CO_\hCS\to
\CO_{\hCS\times_\CX\hCS}(-\ell_\infty,-\ell_\infty+N\ell)[1]\to\ldots
\end{equation}
\eq{conv2}
\ldots\to\CO_{\hCS\times_\CX\hCS}(-\ell_\infty,0)\to\CO_{\hCS\times_\CX\hCS}
\stackrel{\eta}{\longrightarrow}\CO_{\ell_\infty\times_\CX\hCS}
\stackrel{\rho}{\longrightarrow}\CO_{\hCS\times_\CX\hCS}(-\ell_\infty,0)[1]
\to\ldots
\end{equation}
It follows that $\Delta_*\CO_\hCS$ is the convolution of the following
compex of objects of the derived category of coherent sheaves on $\hCS$:
\eq{conv3}
\CO_{\hCS\times_\CX\hCS}(-\ell_\infty,-\ell_\infty+N\ell)
\stackrel{\eta\circ\xi}{\longrightarrow}\CO_{\ell_\infty\times_\CX\hCS}
\stackrel{\rho}{\longrightarrow}\CO_{\hCS\times_\CX\hCS}(-\ell_\infty,0)[1]
\end{equation}
Now since $\CG\simeq\on{pr}_{2*}(\Delta_*\CO_\hCS\otimes\on{pr}_1^*\CG)$, we
see that $\CG$ is the convolution
of the following complex of objects of the derived coherent category of $\hCS$:
\eq{complex}
\pi^*\CG''(-\ell_\infty+N\ell)\stackrel{\alpha}{\longrightarrow}
\pi^*\CG_\infty\stackrel{\pi^*\beta}{\longrightarrow}\pi^*\CG''[1]
\end{equation}
Here $\beta$ enters the exact triangle
\eq{ext}
\ldots\to\CG''\to\CG'\to\imath^*\CG\stackrel{\beta}{\longrightarrow}
\CG''[1]\to\ldots
\end{equation}
while $\alpha$ by adjointness and projection formula is the same
as the direct sum of two morphisms $\alpha':\ \CG''\to\CG_\infty$, and
$\alpha'':\ \CG''(N)\to\CG_\infty$. The condition $\pi^*\beta\circ\alpha=0$
implies $\alpha'=0$.
%while $\alpha''$ is nothing else than $\kappa''_2$ of~\refe{xym}.

Let now $\Delta^\CX$ stand for the diagonal embedding $\CX\hookrightarrow
\CX\times\CX$. Then we have the exact sequences of coherent sheaves on
$\CX\times\CX$:
\eq{reso1}
0\to\CO(0,-N)\oplus\bigoplus_{l=0}^{N-1}\CO(l,-l-1)\to
\bigoplus_{l=0}^N\CO(l,-l)\to\Delta^\CX_*\CO_\CX\to0
\end{equation}
\eq{reso2}
0\to\CO(0,-N-1)\oplus\bigoplus_{l=1}^{N-1}\CO(l,-l-1)\to
\bigoplus_{l=1}^N\CO(l,-l)\to\Delta^\CX_*\CO_\CX\to0
\end{equation}
which yield the resolutions~\refe{comple} and~\refe{compl} for
$\CG''\simeq\on{pr}_{2*}(\Delta^\CX_*\CO_\CX\otimes\on{pr}_1^*\CG'')$.
In particular, $B_l:\ V_l\to V_{l+1}$ is induced by
$y_2:\ H^1(\hCS,\CG(-\ell_\infty+l\ell))\to
H^1(\hCS,\CG(-\ell_\infty+(l+1)\ell))$, and $\delta:\ V_0\to V_N$ is induced
by $y_1:\ H^1(\hCS,\CG(-\ell_\infty))\to H^1(\hCS,\CG(-\ell_\infty+N\ell))$.
%In these terms, $\alpha''$ is nothing else than $\kappa''_1$ of~\refe{xym}.

Furthermore, the morphism $\alpha''(l-N):\ \CG''(l)\to\CG_\infty(l-N)$
induces the morphism $V_l=H^1(\CX,\CG''(l))\to H^1(\CX,\CG_\infty(l-N)=
W_l\oplus\ldots\oplus W_{N-1}$ with components $q_l,\ q_{l+1}B_l,\ldots,\
q_{N-1}B_{N-2}\ldots B_l$. The morphism $\beta(l):\ \CG_\infty(l)\to\CG''(l)[1]$
induces the morphism $W_0\oplus\ldots\oplus W_{l-1}=H^0(\CX,\CG_\infty(l))\to
H^1(\CX,\CG''(l))=V_l$ whose last component is $p_l:\ W_{l-1}\to V_l$.

The exact triangle $\ldots\to\CG''\to\CG'\to\CG_\infty\to\CG''[1]\to\ldots$
along with the acyclicity of $\CG_\infty$ yields an isomorphism
$V_0=H^1(\CX,\CG'')\simeq H^1(\CX,\CG')$. Now $e:\ V_N\to V_0$ is induced
by $\x_2:\ V_N=H^1(\CX,\CG''(N))=H^1(\hCS,\CG(-\ell_\infty+N\ell))\to
H^1(\hCS,\CG)=H^1(\CX,\CG')=V_0$.

Finally, the exact triangle
$\ldots\to\CG''(l-N)\to\CG'(l-N)\to\CG_\infty(l-N)\to\CG''(l-N)[1]\to\ldots$
yields the long exact sequence
\eq{long}
0\to H^1(\CX,\CG''(l-N))\to H^1(\CX,\CG'(l-N))\to H^1(\CX,\CG_\infty(l-N))\to
H^2(\CX,\CG''(l-N))\to\ldots
\end{equation}
A resolution $$0\to\CG''(l-N)\to\CG''\oplus\CG''(l)\oplus\CG''(l+1)\oplus\ldots
\oplus\CG''(N-1)\to\CG''(l+1)\oplus\ldots\oplus\CG''(N)\to0$$ implies
$H^1(\CX,\CG''(l-N))\simeq\on{Ker}\left(V_0\oplus\bigoplus_{m=l}^{N-1}V_m\to
\bigoplus_{m=l+1}^NV_m\right)$ which, together
with~\refe{long}, yields an isomorphism
$$H^1(\CX,\CG'(l-N))\simeq\on{Ker}\left(V_0\oplus\bigoplus_{m=l}^{N-1}V_m
\oplus\bigoplus_{m=l}^{N-1}W_m
\stackrel{\varrho}{\longrightarrow}\bigoplus_{m=l+1}^NV_m\right).$$
Here the ``matrix coefficients'' of $\varrho$ are as follows:
$\delta:\ V_0\to V_N;\ B_m:\ V_m\to V_{m+1};\ \Id:\ V_m\to V_m;\ p_{m+1}:\
W_m\to V_{m+1}$. In particular, we have a morphism
$\varpi:\ H^1(\CX,\CG'(l-N))\to V_l$.

Now the morphism
$x_2((l-N)\ell):\ \CG(-\ell_\infty+l\ell)\to\CG((l-N)\ell)$ gives rise
to the morphism $V_l=H^1(\hCS,\CG(-\ell_\infty+l\ell))\to
H^1(\hCS,\CG((l-N)\ell))=H^1(\CX,\CG'(l-N))$. Composing it with $\varpi$
we obtain the morphism $A_l:\ V_l\to V_l$.

\ssec{ass}{Analysis of stability}
Clearly, at $\ell_\infty$, $C$ is injective, and $D$ is surjective. 
Hence $C$ is injective (as a morphism of sheaves) everywhere, i.e.
$H^{-1}$\refe{mona}$=0$, and the support of $H^1$\refe{mona} does not
intersect $\ell_\infty$. Hence this support lies in the union of the 
exceptional divisor $E$ and finitely many points. Now $H^0$\refe{mona}
cannot have torsion in codimension 2 (as the middle cohomology of a 3-term
complex of vector bundles), so it can only have torsion at curves not
intersecting $\ell_\infty$, i.e. at $E$. If $H^0$\refe{mona} does have
torsion at $E$, then the fibers of $C$ at $E$ (i.e. $x_2=0$) have nontrivial
kernels. We can set $x_1=0$, and then the kernel of $C$ at a point
$[y_1:y_2]\in E$ consists of collections of $v_l\in V_l,\ 0\leq l\leq N$, 
such that $q_lv_l=0,\ 0\leq l\leq N-1;\ ev_N=0;\ B_{N-1}v_{N-1}+y_2v_N=0;\ 
\delta v_0-y_1v_N=0;\ B_{l-1}v_{l-1}+y_2v_l=0=A_lv_l,\ 1\leq l\leq N-1$.
Note that if $v_0=0$ then in case $y_1\ne0$ we obtain $v_N=0$, while in case
$y_2\ne0$ we obtain $v_1=0\ \Rightarrow\ v_2=0\ \Rightarrow\ldots\Rightarrow\
v_N=0$. Moreover, in case $y_2\ne0\ne y_1$ we have
$v_1=y_2^{-1}B_0v_0,\ v_2=y_2^{-2}B_1B_0v_0,\ldots,v_N=y_2^{-N}B_{N-1}\ldots B_1B_0v_0=
y_1^{-1}\delta v_0$. The existence of a nontrivial solution of the equation
$y_2^{-N}B_{N-1}\ldots B_1B_0v_0=y_1^{-1}\delta v_0$ for any $y_1,y_2$ implies that
$\on{Ker}(B_{N-1}\ldots B_1B_0)\cap\on{Ker}(\delta)\ne0$. We take 
$0\ne v_0\in\on{Ker}(B_{N-1}\ldots B_1B_0)\cap\on{Ker}(\delta)$, and 
the corresponding $v_1,\ldots,v_{N-1},v_N=0$. We set $S_0:=\BC v_0\subset V_0,
\ldots,S_{N-1}:=\BC v_{N-1}\subset V_{N-1},\ S_N:=0\subset V_N$. Then 
$S_\bullet\subset V_\bullet$ violates the $\zeta^\bullet$-semistability since
$\dim S_0>\dim S_N$ (see~\cite[Definition~1.1.(1) and~Section~4(ii)]{NaM}).
Conversely, the argument like~\cite[Lemma~7.2.(2)]{NY} proves that  
the fiberwise injectivity of $C$ at the generic point of $E$ implies
$\zeta^\bullet$-semistability. 

The vanishing of $H^1\refe{mona}$ is equivalent to the fiberwise surjectivity
of $D$ everywhere, that is to the fiberwise injectivity of the adjoint
morphism $D^*$ everywhere. Since $D^*$ is clearly injective at $\ell_\infty$
we can study the points off $\ell_\infty$, i.e. we can set $x_1=1$.
The kernel of $D^*$ at a point $(x_2,[y_1:y_2])$ consists of collections
$v_l^*\in V_l^*,\ 0\leq l\leq N$, such that
$p_l^*v_l^*=0,\ 1\leq l\leq N;\ y_1B_0^*v_1^*+y_2\delta^*v_N^*=0;\
(e^*\delta^*-x_2y_1)v_N^*=0;\ e^*B_0^*v_1^*+x_2y_2v_N^*=0;\
B_l^*v_{l+1}^*+y_2v_l^*=0=A_l^*v_l^*-x_2y_1v_l^*,\ 1\leq l\leq N-1$.
Note that $v_0^*$ does not enter these equations, and we can set $v_0=0$.
Given a nontrivial solution of these equations we set 
$S_0^*:=0\subset V_0^*,\ S_1^*:=\BC v_1^*\subset V_1^*,\ldots,S_N^*:=\BC v_N^*
\subset V_N^*$. We set $T_l:=(S_l^*)^\perp\subset V_l,\ 0\leq l\leq N$.
We have $0=\dim S_0^*=\codim T_0\leq\codim T_N=\dim S_N^*$ which violates
the $\zeta^-$-semistability 
(see~\cite[Definition~1.1.(2) and~Section~4(ii)]{NaM}). 
Conversely, the argument like~\cite[Lemma~5.1]{NY} proves that the
fiberwise surjectivity of $D$ implies $\zeta^-$-semistability.

\reft{queer blowup} is proved.

\ssec{l_0 triv}{Parabolic sheaves trivial along $\ell_0$}
We say that a parabolic sheaf $\CF_\bullet$ is {\em trivial along $\ell_0$}
if $\CF_0$ is locally free at $\ell_0$, and $\CF_0|_{\ell_0}$ is trivial; 
moreover, for $-N\leq k\leq 0$ the sheaf $\CF_k$ is locally free at $\ell_0$,
and the quotient sheaves $\CF_k/\CF_{-N},\CF_0/\CF_k$ are both locally free
and trivial (as vector bundles) at $\ell_0$. It is easy to see that the
moduli space $\hCP_{\unl{d},\on{triv}}\subset\hCP_{\unl{d}}$ of parabolic sheaves
trivial along $\ell_0$ is an open subset of $\hCP_{\unl{d}}$, nonempty iff
$d_0=d_1=\ldots=d_N$. 

We also consider an open subset $\hfM_{\unl{d},\on{iso}}\subset\hfM_{\unl{d}}$ 
formed by the representations of $\widehat Q$ such that $B_0,\ldots,B_{N-1}$
are all isomorphisms (evidently, $\hfM_{\unl{d},\on{iso}}$ is nonempty iff 
$d_0=d_1=\ldots=d_N$). The proof of~\reft{queer blowup} admits the following

\cor{naka}
$\Xi(\hfM_{\unl{d},\on{iso}})=\hCP_{\unl{d},\on{triv}}$.
\ecor

\prf
The complex $\CG''\in D(\CX)$ corresponds to the representation~\refe{xquiver}
in $D(K)$. It is easy to check that $B_0,\ldots,B_{N-1}$ in~\refe{xquiver}
are all isomorphisms iff $\Ext_\CX^\bullet(\CG'',\CH)=0$ for any skyscraper sheaf
$\CH$ supported at $0_\CX$. More precisely, $B_l$ is an isomorphism iff
$\Ext_\CX^\bullet(\CG'',\on{Coker}(\CO_\CX(-l-1)\hookrightarrow\CO_\CX(-l)))=0$.
Thus all $B_l$ are isomorphisms iff $\CG''$ has a finite support on $\CX$ disjoint from
$0_\CX$. From the exact triangle~\refe{ext}, $\CG'$ is isomorphic to $\imath^*\CG=\CG_\infty$ 
near $0_\CX$. From the exact triangle~\refe{32}, $\CG$ is isomorphic to $\pi^*\CG_\infty$
near $\ell$, that is $\CF_\bullet$ is trivial along $\ell_0$.
\epr

\ssec{BD}{Blowdown}
The blowdown morphism $\widehat{\BP}{}^2\to\BP^2$ {\em does not} lift to
a morphism of the stacks $\hCS\to\CS'$ (see~\refss{PNN1}). It only gives
rise to a correspondence $\hCS\stackrel{\tilde{\mu}}{\longleftarrow}\CW
\stackrel{\tilde{\nu}}{\longrightarrow}\CS'$.
Since both $\hCS$ and $\CS'$ are toric stacks, this correspondence can be
described in toric terms. Namely, $\hCS$ is given by a fan $\hat{F}$ formed
by the vectors $(1,0);(0,1);(-1,0);(-N,-N)$ in $\BZ^2$, while $\CS'$ is given
by a fan $F'$ formed by the vectors $(1,0);(0,1);(-N,-N)$, and $\CW$ is given
by a fan $\tilde F$ formed by the vectors $(1,0);(0,1);(-N,0);(-N,-N)$.
The evident embedding $F'\subset\tilde{F}$ corresponds to our $\tilde\nu$.
Since $\tilde{F}$ is obtained from $\hat F$ by dilating the vector $(-1,0)$,
we obtain the desired morphism $\tilde\mu:\ \CW\to\hCS$.
We set $\Pi:=\tilde{\nu}_*\tilde{\mu}{}^*:\ D^b\on{Coh}(\hCS)\to
D^b\on{Coh}(\CS')$. Note that the assumptions of~Theorem~4.2(2) of~\cite{Ka}
are satisfied in our situation.

Given $\CF=\CG(\ell)\in\hCP_{\unl{d}}$ the complex $\Pi\CF$ is not
necessarily a torsion free sheaf: it can have the first cohomology (a torsion
sheaf at the origin); it is rather a {\em perverse coherent sheaf}. Its class
is well defined in the Zastava space $\fZ_{\unl{d}}$ (see~\refss{chain}).
Thus we obtain the morphism $\hCP_{\unl{d}}\to\fZ_{\unl{d}}$ which factors
through the morphism $\hfZ_{\unl{d}}\to\fZ_{\unl{d}}$
(since $\fZ_{\unl{d}}$ is affine, $\hfZ_{\unl{d}}$ is normal, and
$\hCP_{\unl{d}}\to\hfZ_{\unl{d}}$ is proper) to be denoted by $\Pi$.

%Comparing the constructions of~\refss{rot} and~\refss{ic} we arrive at the
%following description of $\Pi$.
\conj{PI}
In quiver terms, $\Pi:\ \hfZ_{\unl{d}}\to\fZ_{\unl{d}}$
sends $(V_\bullet,A_\bullet,B_\bullet,\delta,e,p_\bullet,
q_\bullet)$ to $(V'_\bullet,A'_\bullet,B'_\bullet,p'_\bullet,q'_\bullet)$ where
$V'_l:=V_l$ for $l=1,\ldots,N$, and $A'_l:=A_l$ for $l=1,\ldots,N-1$, while
$A'_0:=\delta e$. Furthermore, $B'_l:=B_l$ for $l=1,\ldots,N-1$, while
$B'_0:=B_0e$. Furthermore, $p'_l=p_l$ for $l=1,\ldots,N$, and
$q'_l:=q_l$ for $l=1,\ldots,N-1$, while $q_0:=q_0e$.
\econj

\sec{klein blowup}{Zastava for Kleinian Blowup}

\ssec{KB}{Kleinian Blowup}
We consider $\CS'_1=\BP^2$ with homogeneous coordinates $[z_0:z_1:z_2]$.
We blow up the sheaf of ideals $I$ supported at the origin, and generated
locally by $(z_1^k,z_2)$. This blowup is a singular toric surface lying
in the projectivization $\on{Proj}(\CO(-k)\oplus(\CO(-1))$ of the vector
bundle $\CO(-k)\oplus\CO(-1)$ over $\BP^2$.
%product $\BP^2\times\BP^1$ and defined by equation $z_1^kw=z_2z$ where
%$[z:w]$ are the homogeneous coordinates on $\BP^1$. 
In fact, it has a unique singular point, lying in a chart $\overline{U}{}^2$
with coordinates $z_1,z_2,z$ satisfying
$z_2z=z_1^k$. We define a smooth toric stack $\hCS{}^k_1$ as the stacky
resolution of our blowup at the singular point. The neighbourhood $U^2$ of the
stacky point (the preimage of $\overline{U}{}^2$) is isomorphic to
$\BA^2/\Gamma_k$ (with hyperbolic action). The stack $\hCS{}^k_1$ is given by a fan
$\hat{F}{}^k$ formed by the vectors $(1,0);(0,1);(-1,k-1);(-1,-1)$ in $\BZ^2$.

Let $\ell_0\subset\hCS{}^k_1$ be the proper transform of the line
$\ell_0\subset\CS'_1$ given by the equation $z_2=0$. It lies in the union of
two charts $U^1$ with coordinates $z_1,w$ and $U^\infty$ with coordinates
$z_0,z_2$. We consider the ramified Galois coverings $U^1_N$ with coordinates
$z_1,\sqrt[N]{w}$, and $U^\infty_N$ with coordinates $z_0,\sqrt[N]{z_2}$,
with Galois group $\Gamma_N$. Gluing the stacky quotients $U^1_N/\Gamma_N$,
and $U^\infty_N/\Gamma_N$ with $U^2$, we obtain the smooth toric stack
$\hCS{}^k_N$. It is given by a fan $\hat{F}{}^k_N$ formed by the vectors
$(1,0);(0,1);(-1,k-1);(-N,-N)$ in $\BZ^2$. The preimage of
$\ell_0\subset\hCS{}^k_1$ is denoted by $\ell\subset\hCS{}^k_N$. Its
automorphism group is $\Gamma_N$.

We have a correspondence $\hCS{}^k_N\stackrel{\tilde{\mu}{}_k}{\longleftarrow}
\CW^k\stackrel{\tilde{\nu}{}^k}{\longrightarrow}\CS'_N$ (cf.~\refss{BD})
where $\CW^k$ is given by a fan $\tilde{F}{}^k$ formed by the vectors
$(1,0);(0,1);(-N,N(k-1));(-N,-N)$.
The evident embedding $F'\subset\tilde{F}{}^k$ corresponds to our
$\tilde{\nu}{}^k$.
Since $\tilde{F}{}^k$ is obtained from $\hat{F}{}^k_N$ by dilating the vector
$(-1,k-1)$, we obtain the desired morphism
$\tilde{\mu}{}_k:\ \CW^k\to\hCS{}^k_N$.

\ssec{PSKB}{Parabolic sheaves on Kleinian Blowup}
Let $\hCP{}^k$ be the moduli space of rank $N$ torsion free parabolic sheaves
on $\hCS{}^k_1$ (with parabolic structure along $\ell_0$, and with trivial
first Chern class) trivialized at
infinity; equivalently, $\hCP{}^k$ is the moduli space of rank $N$ torsion
free sheaves on $\hCS{}^k_N$ framed at
$\ell_\infty:\ \imath^*\CF\simeq\CF_\infty:=W_0\otimes\CO_\CX\oplus
W_1\otimes\CO_\CX(-1)\oplus\ldots\oplus W_{N-1}\otimes\CO_\CX(N-1)$
(and with first Chern class trivial off $\ell_\infty$).
Given $\CF\in\hCP{}^k$ the complex
$\Pi^k\CF:=\tilde{\nu}{}^k_*\tilde{\mu}{}^*_k\CF$ is not necessarily a
torsion free sheaf on $\CS'_N$: it can have the first cohomology (a torsion
sheaf at the origin); it is rather a perverse coherent sheaf. Its class is
well defined in the Zastava space $\fZ$ (the direct limit of all
$\fZ_{\unl{d}}$ with respect to the natural embeddings 
$\fZ_{\unl{d}}\hookrightarrow\fZ_{\unl{d}'},\ \unl{d}'\geq\unl{d}$ 
componentwise, adding defect at the origin). Thus we obtain the morphism
$\hCP{}^k\to\fZ$. Our goal is to describe the moduli space $\hCP{}^k$
(in particular, to number its connected components) in quiver terms, as well
as the morphism $\hCP{}^k\to\fZ$. A connected component of $\hCP{}^k$ will
be called {\em good} if it contains a {\em nonempty} open subset formed by
{\em locally free} parabolic sheaves.

\ssec{bicycle}{Rift}
We consider the representations of the following {\em rift} quiver $\widehat{Q}{}^k$:

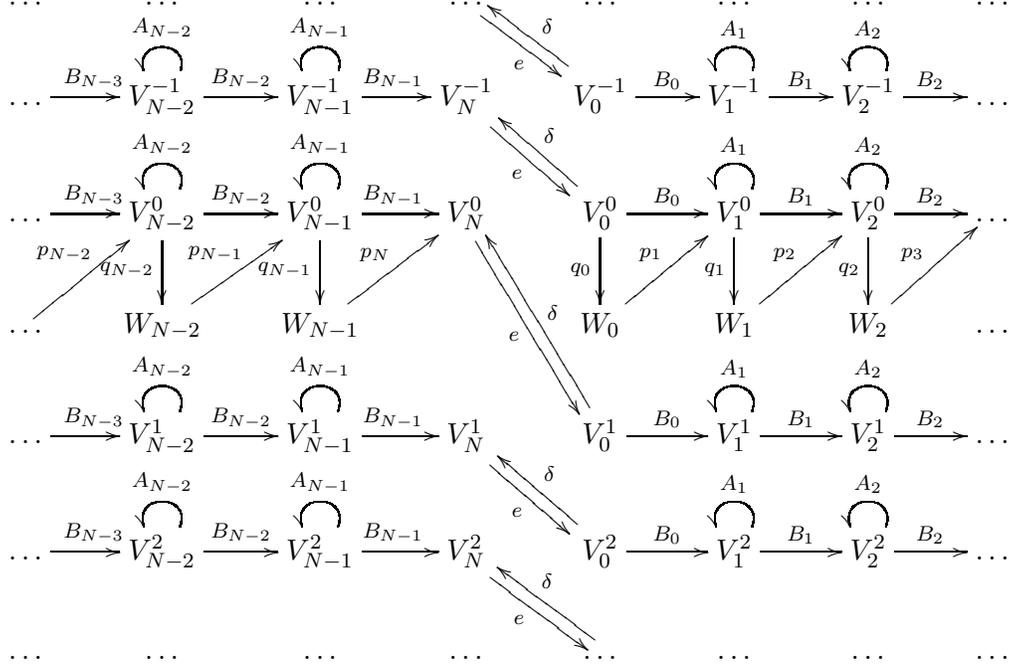
$$\xymatrix{
\ldots & \ldots & \ldots & \ldots \ar@<-.5ex>[dr]_e & \ldots
& \ldots & \ldots & \ldots\\
\ldots \ar[r]^{B_{N-3}}
& V_{N-2}^{-1} \ar@(ur,ul)[]_{A_{N-2}} \ar[r]^{B_{N-2}}
& V_{N-1}^{-1} \ar@(ur,ul)[]_{A_{N-1}} \ar[r]^{B_{N-1}}
& V_N^{-1} \ar@<-.5ex>[dr]_e
& V_0^{-1} \ar@<-.5ex>[ul]_\delta \ar[r]^{B_0}
& V_1^{-1} \ar@(ur,ul)[]_{A_1} \ar[r]^{B_1}
& V_2^{-1} \ar@(ur,ul)[]_{A_2} \ar[r]^{B_2} &\ldots\\
\ldots \ar[r]^{B_{N-3}}
& V_{N-2}^0 \ar@(ur,ul)[]_{A_{N-2}} \ar[r]^{B_{N-2}} \ar[d]_{q_{N-2}}
& V_{N-1}^0 \ar@(ur,ul)[]_{A_{N-1}} \ar[r]^{B_{N-1}} \ar[d]_{q_{N-1}}
& V_N^0 \ar@<-.5ex>[ddr]_e
& V_0^0 \ar@<-.5ex>[ul]_\delta \ar[r]^{B_0} \ar[d]_{q_0}
& V_1^0 \ar@(ur,ul)[]_{A_1} \ar[r]^{B_1} \ar[d]_{q_1}
& V_2^0 \ar@(ur,ul)[]_{A_2} \ar[r]^{B_2} \ar[d]_{q_2}
&\ldots\\
\ldots \ar[ur]^{p_{N-2}}
& W_{N-2} \ar[ur]^{p_{N-1}}
& W_{N-1} \ar[ur]^{p_N}
& & W_0 \ar[ur]^{p_1} & W_1 \ar[ur]^{p_2}
& W_2 \ar[ur]^{p_3}
&\ldots\\
\ldots \ar[r]^{B_{N-3}}
& V_{N-2}^1 \ar@(ur,ul)[]_{A_{N-2}} \ar[r]^{B_{N-2}}
& V_{N-1}^1 \ar@(ur,ul)[]_{A_{N-1}} \ar[r]^{B_{N-1}}
& V_N^1 \ar@<-.5ex>[dr]_e
& V_0^1 \ar@<-.5ex>[uul]_\delta \ar[r]^{B_0}
& V_1^1 \ar@(ur,ul)[]_{A_1} \ar[r]^{B_1}
& V_2^1 \ar@(ur,ul)[]_{A_2} \ar[r]^{B_2} &\ldots\\
\ldots \ar[r]^{B_{N-3}}
& V_{N-2}^2 \ar@(ur,ul)[]_{A_{N-2}} \ar[r]^{B_{N-2}}
& V_{N-1}^2 \ar@(ur,ul)[]_{A_{N-1}} \ar[r]^{B_{N-1}}
& V_N^2 \ar@<-.5ex>[dr]_e
& V_0^2 \ar@<-.5ex>[ul]_\delta \ar[r]^{B_0}
& V_1^2 \ar@(ur,ul)[]_{A_1} \ar[r]^{B_1}
& V_2^2 \ar@(ur,ul)[]_{A_2} \ar[r]^{B_2} &\ldots\\
\ldots & \ldots & \ldots & \ldots & \ldots \ar@<-.5ex>[ul]_\delta
& \ldots & \ldots & \ldots
}$$

Here the upper indices of $V$ run through $\BZ/k\BZ$.
The dimension of $V_l^r$
%depends only on $r$ (i.e. the dimension is constant in a row), and
is denoted by $d_l^r$. We consider the dimension vector
$\hd:=(d_l^r)_{0\leq l\leq N}^{r\in\BZ/k\BZ}$. Furthermore, $\dim W_l=1$,
and all these lines are identified with, say $W_\infty$, so that $W_\infty$ is
the source of all $p_l$ and the target of all $q_l$.

Relations: $0=\delta eB_{N-1}-B_{N-1}A_{N-1}+p_Nq_{N-1}:\ V_{N-1}^0\to V_N^0$.

$0=A_1B_0-B_0e\delta+p_1q_0:\ V_0^0\to V_1^0$.

$0=A_{l+1}B_l-B_lA_l+p_{l+1}q_l:\ V_l^0\to V_{l+1}^0$ for $l=1,\ldots,N-2$.

$0=\delta eB_{N-1}-B_{N-1}A_{N-1}:\ V_{N-1}^r\to V_N^r$ for $r\ne0$.

$0=A_1B_0-B_0e\delta:\ V_0^r\to V_1^r$ for $r\ne0$.

$0=A_{l+1}B_l-B_lA_l:\ V_l^r\to V_{l+1}^r$ for $l=1,\ldots,N-2,\ r\ne0$.

We denote a typical representation of $\widehat{Q}{}^k$ by $Y$.
We denote by $\sM_\hd$ the scheme of representations of $\widehat{Q}{}^k$ of
dimension $\hd$. We denote by $G_\hd$ the group
$\prod_{0\leq l\leq N}^{r\in\BZ/k\BZ}\on{GL}(V_l^r)$; it acts naturally on
$\sM_\hd$.

\ssec{sc}{Stability conditions}
Following~\cite[Section 4(ii)]{NaM} we consider the enhanced dimension vector
$\td:=(d_l^r,1)$ with one extra coordinate equal to $\dim W_\infty=1$.
We consider a vector $\zeta^\bullet=(\zeta_l^r)$
where $\zeta_N^r=-1,\ \zeta_0^r=1,\ \zeta_l^r=0$ for $1\leq l\leq N-1$.
Also, for $0<\varepsilon\ll1$ we consider $\zeta^-:=\zeta^\bullet-(\varepsilon,
\ldots,\varepsilon)$. We set $\zeta^-_\infty:=-\langle\zeta^-,\hd\rangle$,
and $\zeta^\bullet_\infty:=-\langle\zeta^\bullet,\hd\rangle$ where
$\langle\cdot,\cdot\rangle$ stands for the sum of products of coordinates
(the standard scalar product). Finally, we set
$\tzeta^-:=(\zeta^-,\zeta^-_\infty)$, and
$\tzeta^\bullet:=(\zeta^\bullet,\zeta^\bullet_\infty)$.

For a nonzero $\widehat{Q}{}^k$-submodule $Y'\subset Y$ of enhanced dimension $\td'$
(where the last coordinate may be either 1 or 0) we define the slope by
$$\theta^-(Y'):=\frac{\langle\tzeta^-,\td'\rangle}
{\langle(1,\ldots,1),\td'\rangle},\
\theta^\bullet(Y'):=\frac{\langle\tzeta^\bullet,\td'\rangle}
{\langle(1,\ldots,1),\td'\rangle}.$$
We say that a $\widehat{Q}{}^k$-module $Y$ is $\zeta^-$-semistable (resp.
$\zeta^\bullet$-semistable) if for any nonzero submodule $Y'\subset Y$
we have $\theta^-(Y')\leq\theta^-(Y)$ (resp.
$\theta^\bullet(Y')\leq\theta^\bullet(Y)$). We say $Y$ is
$\zeta^-$-stable (resp. $\zeta^\bullet$-stable) if the inequality is strict
unless $Y'=Y$. Note that $\zeta^-$-stability is equivalent to
$\zeta^-$-semistability.

We define a scheme $\hfM_\hd^k$ as the moduli space
of $\zeta^-$-semistable (equivalently, $\zeta^-$-stable) $\widehat{Q}{}^k$-modules.
By GIT, $\hfM_\hd^k$ is the projective spectrum of the ring of
$\hG_\hd$-semiinvariants in $\BC[\hsM_\hd]$.
Furthermore, we define a scheme $\hfZ_\hd^k$ as the moduli space of
$S$-equivalence classes of $\zeta^\bullet$-semistable $\widehat{Q}^k$-modules.
Since the stability condition $\zeta^\bullet$ lies on a wall of the chamber
containing $\zeta^-$, we have a projective morphism $\pi_{\zeta^\bullet,\zeta^-}:\
\hfM_\hd^k\to\hfZ_\hd^k$.

\th{maink}
A good connected component of $\hCP^k$ is isomorphic to $\hfM_\hd$ for
a dimension vector $\hd$ such that $d_0^0=d_N^0$, and for
$r\ne0,\ d_l^r=d_m^r\ \forall l,m=1,\ldots,N$.
\eth

The proof is given in the next Subsection.

\ssec{Gammak}{The action of $\Gamma_k$}
The action of $\Gamma_k$ on $\CS'_1$ (see~\refss{gammak}) lifts to the
action of $\Gamma_k$ on the blowup $\hCS_1$, and also lifts to the action
of $\Gamma_k$ on $\hCS_N$. Hence $\Gamma_k$ acts on the moduli space
$\hCP_{\unl{d}}$ of parabolic sheaves on $\hCS_1$ trivialized at infinity via
its action on $\hCS_1$ and the {\em trivial} action on the trivialization at
infinity.

The fixed point variety
$\hCP_{\unl{d}}^{\Gamma_k}=\hfM_{\unl{d}}^{\Gamma_k}$ can be described in
quiver terms as well. Namely, the construction of quiver in~\refss{ic} implies
that the action of the generator $\zeta_k$ of $\Gamma_k$ on the quiver
components works as follows:
$\zeta_{kN}(A_\bullet,B_\bullet,e,\delta,p_\bullet,q_\bullet)=
(A_\bullet,\zeta_kB_\bullet,\zeta_k^{-2}e,\zeta_k^2\delta,\zeta_kp_\bullet,
q_\bullet)$. It follows that the various connected components of
$\hfM_{\unl{d}}^{\Gamma_k}$ are isomorphic to $\hfM_\hd^k$ for various
dimension vectors $\hd$ such that $\unl{d}=(\sum_{r\in\BZ/k\BZ}d_0^r,
\sum_{r\in\BZ/k\BZ}d_1^r,\ldots,\sum_{r\in\BZ/k\BZ}d_N^r)$.

Among these connected components we single out the ones classifying the
parabolic sheaves with $\Gamma_k$-equivariantly trivial determinant. This
is the condition $d_0^r=d_N^r$ for any $r\in\BZ/k\BZ$.
Also, we single out the components classifying the $\Gamma_k$-equivariant
parabolic sheaves with the {\em trivial} defect class
$[\overline{\on{def}}(\CF_\bullet)]$ (see~\refss{def}).
This is the condition that for
$r\ne0,\ d_l^r=d_m^r\ \forall l,m=1,\ldots,N$. We will refer to such connected
components (satisfying both of the above conditions) as the {\em admissible}
ones.

Now let us consider the stacks $\hCS_1/\Gamma_k$, and $\hCS_1/|\Gamma_k$.
The latter stands for the coarse (categorical) quotient in a neighbourhood
of $\ell_0$ (but the stacky quotient elsewhere). We have an evident projection
$\phi:\ \hCS_1/\Gamma_k\to\hCS_1/|\Gamma_k$.
Given a $\Gamma_k$-equivariant parabolic sheaf $\CF$ on $\hCS_1$ lying
in an admissible connected component $\hfM_\hd^k$, the parabolic sheaf
$\phi_*\CF$ is trivialized at infinity. Similarly to~\refss{dir},
$\phi_*$ induces an isomorphism of an admissible connected component
$\hfM_\hd^k$ with a connected component of the moduli space of torsion free
parabolic sheaves on $\hCS_1/|\Gamma_k$ (with the trivial action of $\Gamma_k$
at the trivialization at infinity). However, the stacks $\hCS_1/|\Gamma_k$
and $\hCS_1^k$ are isomorphic off infinity, so the latter connected component
is nothing else than a connected component of $\hCP^k$.

This completes the proof of~\reft{maink}.

\ssec{pik}{The morphism $\Pi^k$}
For an admissible dimension vector $\hd$ we denote the corresponding connected
component of $\hCP^k$ by $\hCP^k_\hd\simeq\hfM_\hd^k$. We are going to describe
in quiver terms the morphism $\hCP^k_\hd\to\fZ$ of~\refss{PSKB}.
Note that since $\fZ$ is affine, $\pi_{\zeta^\bullet,\zeta^-}:\
\hfM_\hd^k\to\hfZ_\hd^k$ is proper, and $\hfZ_\hd^k$ is normal, the morphism
$\hCP^k_\hd\to\fZ$ factors as the composition of $\pi_{\zeta^\bullet,\zeta^-}$ and
a certain morphism $\hfZ_\hd^k\to\fZ$ to be denoted by $\Pi^k$.

The comparison of constructions of~\refss{Gammak},~\refss{dir},~\refss{BD},
and~\refco{PI}
implies that $\Pi^k=\psi^k\circ\Psi^k$ where
$\Psi^k:\ \hfZ^k_\hd\to\fZ_{\widetilde{\unl{d}}}$ is defined as follows.
First, $\widetilde{\unl{d}}:=(d_l^r)_{1\leq l\leq N}^{r\in\BZ/k\BZ}$ is obtained
from the vector $\hd$ just by erasing the coordinates $d_0^r,\ r\in\BZ/k\BZ$.
Second, $\Psi^k$ acts on the quiver data as follows: $'V_l^r:=V_l^r$ for
$1\leq l\leq N,\ r\in\BZ/k\BZ$. Furthermore, $'B_N:=B_0e:\ V_N^r\to V_1^{r+1}$
and $'A_N:=\delta e:\ V_N^r\to V_N^r$ for $r\in\BZ/k\BZ$. Furthermore,
$'q_0:=q_0e:\ V_N^{-1}\to W_0$, and all the other primed letters are equal
to the corresponding letters without primes.

In particular, we see that $\Pi^k(\hfZ^k_\hd)$ lands into the connected
component $\fZ_{\unl{d}}$ where $\unl{d}=(d_1^0,\ldots,d_N^0)$.

\ssec{open}{An open piece}
We consider the following (admissible) dimension vector: $d_l^r=v_r$ for any
$l=0,\ldots,N$. Let $\hfM_{\hd,\on{iso}}^k\subset\hfM_\hd^k$ be the open subset given
by the condition that $B_{N-1}B_{N-2}\ldots B_1B_0:\ V_0^r\to V_N^r$ is
an isomorphism for any $r\in\BZ/k\BZ$ (equivalently, all the $B_l$ are
isomorphisms). It follows from~\refc{naka} that
this open subset classifies the parabolic sheaves on
$\hCS^k_1$ {\em trivial along} $\ell_0$. Then the trivialization at infinity
extends through $\ell_0$ as well, and we are left with a torsion free sheaf
on the open set $U^2$ (notations of~\refss{KB}) trivialized at infinity.
According to~\cite[4.2]{Nak}, the moduli space of torsion free sheaves on $U^2$
trivialized at infinity is the classical Nakajima quiver variety $\fM(v,w)$ of
type $\widetilde{A}_{k-1}$ where $w=(N,0,\ldots,0)$, and
$v=(v_0,\ldots,v_{k-1})$.

The isomorphism $\Phi:\ \hfM_{\hd,\on{iso}}^k\iso\fM(v,w)$
in quiver terms is given by $W''_0:=W_0\oplus\ldots\oplus W_{N-1}$, and
$W''_r:=0$ for $r\ne0$. Furthermore, $V''_r:=V_0^r$.
Furhermore, $B''_r:=eB_{N-1}B_{N-2}\ldots B_1B_0:\ V''_r\to V''_{r+1}$,
while $A''_r:=(B_{N-1}B_{N-2}\ldots B_1B_0)^{-1}\delta:\ V''_r\to V''_{r-1}$.
Finally, $p''_0:=\bigoplus_{1\leq l\leq N}p^l:\
W''_0\to V''_0$, and $q''_0:=\bigoplus_{0\leq l\leq N-1}q^l:\ V''_0\to W''_0$,
where
$p^l:=B_0^{-1}B_1^{-1}\ldots B_{l-2}^{-1}B_{l-1}^{-1}p_l:\ W_l\to V_0^0$, and
$q^l:=q_lB_{l-1}B_{l-2}\ldots B_1B_0:\ V_0^0\to W_l$.

\ssec{numer}{The image of $\Pi^k$}
Note that while~\reft{maink} provides the necessary (admissibility)
conditions for a component $\hCP{}^k_\hd$ to be good, it does not give the
sufficient conditions. Let us give such sufficient conditions in the setup
of~\refss{open}.
Thus we restrict ourselves to the components $\hfM_\hd^k$ which contain the
{\em nonempty} open subset formed by the {\em locally free} parabolic sheaves.
Equivalently, we are interested in the components $\hfM_{\hd,\on{iso}}^k\simeq\fM(v,w)$
which contain the {\em nonempty} open subset formed by the vector bundles
on $U^2$. The corresponding dimension vectors $(d_l^r)$ (equivalently,
vectors $v$), will be called {\em good}.
The well-known Nakajima criterion states that $v$ is good iff
the $\widehat{\mathfrak{sl}(k)}$-weight $w-Cv$ is dominant and has nonzero
multiplicity in the level $N$ vacuum integrable module $L(w)$. Here $C$ is the
affine Cartan matrix of $\tilde A_{k-1}$. More explicitly,
the dominance condition reads as follows:
$v_0+v_2\geq2v_1,\ldots,v_{k-2}+v_0\geq2v_{k-1},v_{k-1}+v_1+N\geq2v_0$.
In particular, $v_0\geq v_i\ \forall i\in\BZ/k\BZ$.

Let us note that the type of our bundles on $U^2$ at the hyperbolic point is
$\blambda=\ ^t(w-Cv)$ where the transposition
is extensively discussed in~\cite{BF},~\cite{BF2}. In notations of
{\em loc. cit.} (especially~Section~7 of~\cite{BF}), we have
$\fM^{\on{reg}}(v,w)\simeq\on{Bun}^\lambda_{\SL(N),\mu}(\BA^2/\Gamma_k)$
where $\lambda,\mu$ are the following integrable
$\mathfrak{sl}(N)_{\on{aff}}$-weights: $\lambda=(k,\blambda,
\frac{a+\frac{(\bmu,\bmu)}{2}-\frac{(\blambda,\blambda)}{2}}{k}),\
\mu=(k,0,0)$. Here $a$ stands for the second Chern class of the
$\Gamma_k$-equivariant vector bundles on a compactification of $\BA^2$.
According to~\cite{BF}, the Nakajima criterion can be equivalently reformulated
as follows: $v$ is good iff $\mu$ has nonzero multiplicity in the
level $k$ integrable $\mathfrak{sl}(N)_{\on{aff}}$-module $L(\lambda)$.
In particular, $\lambda\geq\mu$, i.e. the difference
$\alpha:=\lambda-\mu$ is a linear combination of simple roots of
$\mathfrak{sl}(N)_{\on{aff}}$ with coefficients in $\BN$.

We view $\alpha$ as a vector with coordinates $(\alpha_1,\ldots,\alpha_N)$.
It is easy to see that $(\alpha_1,\ldots,\alpha_N)\leq(v_0,\ldots,v_0)$
componentwise, and hence we have an embedding $\fZ_\alpha\hookrightarrow
\fZ_{(v_0,\ldots,v_0)}$ adding the defect of the complementary degree at the
origin.

From now on let us write $\hfZ^\lambda_\mu$ for $\hfZ_\hd^k$.

\conj{psycho}
Consider a good component $\hfZ^\lambda_\mu=\hfZ_\hd^k
\stackrel{\Pi^k}{\longrightarrow}\fZ_{v_0,\ldots,v_0}$.

(a) The image of $\Pi^k$ is contained in
$\fZ_\alpha\subset\fZ_{v_0,\ldots,v_0}$, so we may and will view $\Pi^k$
as a morphism $\hfZ^\lambda_\mu\to\fZ_\alpha$.

(b) The morphism $\Pi^k:\ \hfZ^\lambda_\mu\to\fZ_\alpha$ is birational and
stratified semismall, so that the direct image $\Pi^k_*\on{IC}(\hfZ^\lambda_\mu)$
is a direct sum of IC-sheaves of certain strata of $\fZ_\alpha$ with
certain multiplicities.

(c) For $\beta\leq\alpha$, and the corresponding stratum
$\fZ_\beta\subset\fZ_\alpha$, the multiplicity $m_\beta$ of
$\on{IC}(\fZ_\beta)$ in $\Pi^k_*\on{IC}(\hfZ^\lambda_\mu)$ equals
the weight multiplicity $L^\lambda(\lambda-\alpha+\beta)$ of the
integrable $\mathfrak{sl}(N)_{\on{aff}}$-module $L^\lambda$.
\econj

\sec{general G}{General $G$}

\ssec{closure}{Arbitrary groups} Let $G$ be an almost simple simply connected
group with the Lie algebra $\fg$. We have the adjoint representation
$G\to\SL(\fg)=\SL(N)$ where $N=\dim\fg$. We choose a Borel subgoup $B_N$ of
$\SL(N)$ containing the image of the positive Borel subgoup $B\subset G$.

We consider the moduli space of $G$-bundles on $\hCS_1/|\Gamma_k$
(see~\refss{Gammak}) equipped with a reduction to $B$ along $\ell_0$ and
framing at $\ell_\infty$. The component of this moduli space having an open
piece $\Bun^\lambda_{G,\mu}(\BA^2/\Gamma_k)$ (cf.~\refss{numer}) will be
denoted by $\widetilde{Z}{}^\lambda_{G,\mu}$. Here $\lambda,\mu\in
\Lambda^+_{\on{aff},k},\ \lambda\geq\mu=(k,0,0)$. The adjoint homomorphism
$\on{ad}:\ (G,B)\to(\SL(N),B_N)$ induces a closed embedding
$\on{ad}:\ \widetilde{Z}{}^\lambda_{G,\mu}\hookrightarrow
\widetilde{Z}{}^{\on{ad}_*\lambda}_{\SL(N),\on{ad}_*\mu}$. We define
$\overline{\CG Z}{}^{\lambda,\alpha}_{G_{\aff}}$ as the closure of
$\on{ad}(\widetilde{Z}{}^\lambda_{G,\mu})$ in
$\hfZ^{\on{ad}_*\lambda}_{\on{ad}_*\mu}\supset\widetilde{Z}{}^{\on{ad}_*\lambda}_{\SL(N),\on{ad}_*\mu}$; here $\alpha:=\lambda-\mu$. Then $\Pi^k$ restricts to the proper morphism
$\phi:\ \overline{\CG Z}{}^{\lambda,\alpha}_{G_{\aff}}\to Z^\alpha_{G_{\aff}}$.

\conj{psi chi}
(a) The morphism $\phi:\ \overline{\CG Z}{}^{\lambda,\alpha}_{G_{\aff}}\to Z^\alpha_{G_{\aff}}$ is birational
and stratified semismall, so that the direct image
$\phi_*\on{IC}(\overline{\CG Z}{}^{\lambda,\alpha}_{G_{\aff}})$
is a direct sum of IC-sheaves of certain strata of $Z^\alpha_{G_{\aff}}$ with
certain multiplicities.

(b) For $\beta\leq\alpha$, and the corresponding stratum
$Z^\beta_{G_{\aff}}\subset Z^\alpha_{G_{\aff}}$, the multiplicity $m_\beta$ of
$\on{IC}(Z^\beta_{G_{\aff}})$ in $\phi_*\on{IC}(\overline{\CG Z}{}^{\lambda,\alpha}_{G_{\aff}})$ equals
the weight multiplicity $L^\lambda(\lambda-\alpha+\beta)$ of the
integrable $G^\vee_{\on{aff}}$-module $L^\lambda$.
\econj

\ssec{rep}{Repellents}
As we have mentioned in the Introduction,~\refco{psi chi} is an affine analogue of the
(known) statement about the convolution between the (classical) affine Grassmannian and
Zastava of $G$. This classical statement is deduced from the description of fibers of
the convolution morphism as the intersections of Schubert varieties in the affine 
Grassmannian with {\em semiinfinite orbits}. We expect a similar description applies in
the affine situation. We will define an open subset of 
$\overline{\CG Z}{}^{\lambda,\alpha}_{G_{\aff}}$ which identifies with a transversal slice
in the double affine Grassmannian of $G$, and prove that its intersection with an
affine analogue of a semiinfinite orbit lies in the central fiber of $\phi$.

We have an intermediate open subset $\widetilde{Z}{}^{\on{ad}_*\lambda}_{\SL(N),\on{ad}_*\mu}\subset
\check{Z}{}^{\on{ad}_*\lambda}_{\SL(N),\on{ad}_*\mu}\subset\hfZ^{\on{ad}_*\lambda}_{\on{ad}_*\mu}$ specified in quiver terms by the condition that the composition $B_{N-1}B_{N-2}\ldots B_1B_0:\ V_0^r\to V_N^r$ is an
isomorphism for any $r\in\BZ/k\BZ$, cf.~\refss{open}. It is nothing else than the Uhlenbeck space
$\CU^{\on{ad}_*\lambda}_{\SL(N),\on{ad}_*\mu}(\BA^2/\Gamma_k)$. The closure of
$\on{ad}(\widetilde{Z}{}^\lambda_{G,\mu})$ in $\CU^{\on{ad}_*\lambda}_{\SL(N),\on{ad}_*\mu}(\BA^2/\Gamma_k)=
\check{Z}{}^{\on{ad}_*\lambda}_{\SL(N),\on{ad}_*\mu}\supset\widetilde{Z}{}^{\on{ad}_*\lambda}_{\SL(N),\on{ad}_*\mu}$ is nothing else than the Uhlenbeck space $\CU^\lambda_{G,\mu}(\BA^2/\Gamma_k)=\ol\CW{}^\lambda_{G_{\aff},\mu}$.
In~Section~3.2 of~\cite{BF3} we have introduced the locally closed subvariety $\fT^e_\mu\subset\CU^\lambda_{G,\mu}(\BA^2/\Gamma_k)=\ol\CW{}^\lambda_{G_{\aff},\mu}$ as the repellent of a certain $\BC^*$-action ($e$ stands for the neutral element of the affine Weyl group). 
(In type $A$ these repellents were introduced in~\cite{Na} under the name of 
{\em MV cycles}.) We conjecture that
the central fiber $\phi^{-1}(i^\alpha_\alpha(0))\cap\CU^\lambda_{G,\mu}(\BA^2/\Gamma_k)$ coincides with
the repellent $\fT^e_\mu\subset\CU^\lambda_{G,\mu}(\BA^2/\Gamma_k)$. We can prove only one inclusion.

First, we recall the definition of the $\BC^*$-action for the reader's 
convenience. We choose a Cartan torus $T\subset B\subset G$ with the coweight
lattice $\Lambda$. The 2-dimensional torus $\BC^*\times\BC^*$ acts on $\BA^2$
naturally: $(a,b)\cdot(z,t)=(az,bt)$, and on $\CU^\lambda_{G,\mu}(\BA^2/\Gamma_k)$
by the transport of structure. 
Let $\BC^*_{\on{hyp}}:=\{(c,c^{-1})\}\subset\BC^*\times\BC^*$ 
%(resp.  $\BC^*_\Delta:=\{(c,c)\}\subset\BC^*\times\BC^*$) 
stand for the antidiagonal, alias hyperbolic, 
%(resp. diagonal) 
subgroup. Let us denote the torus $\BC^*_{\on{hyp}}\times T$ by $\widehat T$. 
This is a Cartan torus of $\widetilde G$ with the coweight lattice 
$\widehat\Lambda$.
Thus, the slices $\ol\CW{}^\lambda_\mu=\CU^\lambda_{G,\mu}(\BA^2/\Gamma_k)$ are equipped with the action of the torus $\widehat{T}$.
Let $I$ (resp. $I_{\aff}=I\sqcup i_0$) stand for the set of vertices of the Dynkin diagram of $G$ (resp. of $G_{\aff}$). For $i\in I$ we denote by 
$\ol\omega_i\in\Lambda$ the corresponding fundamental coweight of $G$,
and we denote by $a_i\in\BN$ the corresponding label of the Dynkin diagram of $G_{\aff}$.
Then $\omega_{i_0}:=(1,0)\in\BZ\times\Lambda=\widehat\Lambda,\ \omega_i:=(a_i,\ol\omega_i)\in\widehat\Lambda,\ i\in I$, are the fundamental coweights of $G_{\aff}$. We set $\rho:=\sum_{i\in I_{\aff}}\omega_i$ (not to be confused with the halfsum $\bar\rho$ of positive coroots of $G$). 
%The group $W_{\aff,k}$ acts on 
%$\widehat\Lambda$, and for $w\in W_{\aff,k}$ let us view $w\rho$ as a 
%one-parametric subgroup $\BC^*\to\widehat T$. 
The torus $\widehat T$ acts on 
$\ol\CW{}^\lambda_\mu$ with the only fixed point, to be denoted abusively by $\mu$, and we define $\fT^e_\mu\subset\ol\CW{}^\lambda_\mu$ as the repellent
$\fT^e_\mu:=\{g\in\ol\CW{}^\lambda_\mu:\ \lim_{c\to\infty}\rho(c)g=\mu\}$.

\prop{repel}
$\fT^e_\mu\subset\phi^{-1}(i^\alpha_\alpha(0))\cap\CU^\lambda_{G,\mu}(\BA^2/\Gamma_k)$.
\eprop

\prf
If a point $a$ lies in $\fT^e_\mu$, then $\phi(a)$ is repelled from $i^\alpha_\alpha(0)\in Z^\alpha_{G_{\aff}}$ under
the following action of $\BC^*$ on $Z^\alpha_{G_{\aff}}$. Recall that $Z^\alpha_{G_{\aff}}$ is a certain closure of the moduli space of $G$-bundles on $\BP^2$ (with homogeneous coordinates $[z_0:z_1:z_2]$) trivialized at $\ell_\infty$ (given by $z_0=0$) and equipped with a reduction to $B$ along $\ell_0$ (given by $z_2=0$). The Cartan torus $T\subset B\subset G$ acts on $Z^\alpha_{G_{\aff}}$ via trivialization at $\ell_\infty$, while $\BC^*_{\on{vert}}$ acts on
$\BP^2$ by $c[z_0:z_1:z_2]=[z_0:z_1:cz_2]$, and hence on $Z^\alpha_{G_{\aff}}$ by transport of structure.
Note that the action of $\BC^*_{\on{vert}}$ lifts to the action on 
$\widehat\BP{}^2_k$ with the following property: if $f$ is a non $\BC^*_{\on{vert}}$-fixed point on the exceptional divisor, then as $c\in\BC^*_{\on{vert}}$ tends to infinity, $c\cdot f$ tends to the singular point of the exceptional divisor. 
Moreover, this action of $\BC^*_{\on{vert}}$ on $\widehat\BP{}^2_k$ in the chart
$\overline{U}{}^2$ (notations of~\refss{KB}) coincides with the action of
$\BC^*_{\on{hyp}}$ on $\BA^2/\!/\Gamma_k$.
We consider the one-parametric subgroup $\BC^*\to\BC^*_{\on{vert}}\times T:\ 
c\mapsto(c^h,\bar\rho(c))$ where
$\bar\rho$ is the halfsum of positive coroots of $G$ viewed as a cocharacter of $T$, while $h$ is the Coxeter number of $G$. The desired action of $\BC^*$ on $Z^\alpha_{G_{\aff}}$ is the action of this one-parametric subgroup. 

Finally, the only points of $Z^\alpha_{G_{\aff}}$ repelled from anything 
at all are the $T$-fixed points $\BA^\alpha\subset Z^\alpha_{G_{\aff}}$. 
In effect, we have the projection $\varrho: Z^\alpha_{G_{\aff}}\to\BA^\alpha$
(see~\cite[Section~9]{BFG}) equivariant with respect to the $\BC^*$-action.
According to the factorization principle~\cite[Corollary~9.4]{BFG}, it suffices
to prove the desired statement for the points of $Z^\alpha_{G_{\aff}}$ lying in
the central fiber $\overline\CF{}^\alpha:=\varrho^{-1}(\alpha\cdot0)$. This
is clear from the stratification~\cite[(7) in~Section~4.1 (or (4.2) 
in~Section~4.1 of arXiv:0912.5132)]{BFK} of 
$\overline\CF{}^\alpha$.

It follows $\phi(a)=i^\alpha_\alpha(0)$.
\epr

\bigskip
\footnotesize{
{\bf A.B.}: Department of Mathematics, Brown University,
151 Thayer St., Providence RI
02912, USA;\\
{\tt braval@math.brown.edu}}

\footnotesize{
{\bf M.F.}: IMU, IITP and National Research University Higher School of 
Economics\\
Department of Mathematics, 20 Myasnitskaya st, Moscow 101000 Russia;\\
{\tt fnklberg@gmail.com}}

%\footnotesize{
%{\bf A.K.}: Algebra Section, Steklov Mathematical Institute, 8 Gubkin st,
%Moscow 119991 Russia;\\
%The Poncelet Laboratory, Independent University of Moscow;\\
%{\tt akuznet@mi.ras.ru}}

\end{document}